\documentclass[12pt]{article}
\usepackage{geometry}
\usepackage{graphicx}
\usepackage{verbatim}
\usepackage{array}
\usepackage{xcolor}
\usepackage{epstopdf}
\usepackage{amsmath,amsfonts,amssymb}
\usepackage{amsthm}
\usepackage{natbib}
\usepackage{multirow}
\usepackage{setspace}
\usepackage{tabu}
\usepackage{authblk}
\usepackage{caption}
\usepackage{graphicx}
\newlength{\myMheight}

\usepackage{multicol}
\usepackage[font=scriptsize]{subfig}
\newenvironment{keywords}{%ffor
\begin{changemargin}{1cm}{1cm}
\noindent{\bf Keywords:}
}
{\end{changemargin} }

\newtheoremstyle{teorema}{\topsep}{\topsep}{}{}{}{ }{\newline}
{\fcolorbox{white}{white}{\bf \color{black}
\textsl{\thmname{#1} \thmnumber{#2}} \thmnote{\{#3\}}}}

\theoremstyle{teorema}

\begin{document}
\onehalfspacing

%%%%%%%%%%%%%%%%%%%%
%%%%%%% Front page %%%%%%%
%%%%%%%%%%%%%%%%%%%%
%%\begin{frontmatter}
\title{A new Ordinal Regression procedure for Multiple Criteria Decision Aiding: the case of the space time model for a sustainable Ecovillage}

\author[1]{Maria Barbati}
\author[2]{Salvatore Greco}
\author[3]{Isabella M. Lami}

\affil[1]{Department of Economic, Ca' Foscari University of Venice, Venice, Italy}
\affil[2]{Department of Economics and Business, University of Catania,   Catania, Italy}
\affil[3]{Interuniversity Department of Regional and Urban Studies and Planning, Politecnico di Torino, Turin, Italy}

\date{}
\maketitle

%\vspace{-1cm}

%%%%%%%%%%%%%%%%%%%%(
\begin{abstract}
In this paper, we present a methodology for handling territorial planning decisions. The proposed approach is based on a multiobjective optimization suggesting which facility to implement, in which location, and at which time.  In  this context, we define a new elicitation procedure to handle Decision Makers (DMs) preferences with an intrinsic and more general interest that goes beyond the specific decision problem. In particular, the user's preferences are elicited by conjugating the deck of cards method with the ordinal regression approach allowing the DM to provide preference information in terms of ranking and pairwise comparing with regard to the intensity of preference of some solutions of the optimization problem. Then, the score of the reference solutions obtained through the deck of the cards method is used as a basis for an ordinal regression procedure that, to take into account interaction between criteria, represents DM's multicriteria preferences by means of a value function expressed in terms of a Choquet Integral. The obtained value function is then used to define a multiobjective optimization problem. The new feasible solutions obtained by the resolution of the optimization problem are proposed to the DM to verify his appreciation and collect new preference information to iterate the interaction procedure ending when the DM is satisfied of the proposed solution.
We apply our methodology to a real world problem to handle the planning procedure of a sustainable Ecovillage in the province of Turin (Italy). We consider a set of facilities to be distributed in a given space in a proper temporal sequence that we conveniently formulated in terms of the space-time model introduced in \cite{barbati2020general}. We interact with the President of the cooperative owning the Ecovillage to detail what  facilities  of the Ecovillage should be selected among the proposed ones, where they should be located and when they should be planned.

\end{abstract}

%evaluated with respect to a multiplicity of criteria.

\providecommand{\keywords}[1]
{
  \small	
  \textbf{\textit{Keywords---}} #1
}

\begin{Keywords}
Territorial Planning; Ordinal Regression; Deck of the Cards Method; Multiobjective Optimization
\end{Keywords}

%%\end{frontmatter}

\pagenumbering{arabic}
%%\journal{Omega}

\section{Introduction}
Decisions regarding territorial planning are very complex for the many points of view that must be taken into account and for the many actors to be involved. In this perspective, transparent and participative procedures would be beneficial for decision support in this domain. With this aim, we propose a Multiple Criteria Decision Aiding (MCDA) methodology (for a collection of state-of-the-art surveys, see \cite{greco2016multiple}) which formulates this type of problem in terms of a multiobjective optimization problem in which specific attention is paid to the aspects related to the interaction with the Decision Maker (DM). In particular, we propose a preference elicitation procedure that conjugates two main approaches in this domain, the deck of the cards method \citep{figueira2002determining,abastante2020introduction} and the ordinal regression \citep{jacquet1982assessing,jacquet2001preference}. The deck of the cards method permits the DM to express his preferences in a straightforward and understandable form, while the ordinal regression permits to induce parameters of the adopted decision model effectively. With respect to the basic model of ordinal regression, the advantage of the proposed methodology is the consideration not only of ordinal information of the type ``alternative $a$ is preferred to alternative $b$'', but also of more cardinal information of the type  ``$a$ is more preferred to $b$, than $c$ is preferred to $d$'', that, just thanks to the deck of cards method, can be handled with a very ``user-friendly'' procedure. We shall call this new methodology Deck of the cards-based Ordinal Regression (DOR).

 Let us explain in more detail how we articulate these ideas in the proposed methodology. We approach territorial planning in terms of multiobjective optimization (Miettinen et al., 2008). In fact, planning the city and the territory generally means having to make decisions regarding the choice of facilities, their location, and their timing under some constraints (Pujadas et al., 2017). For this purpose, several objectives must be considered, and, in general, there is no best solution for all the objectives considered. In this context, a solution is considered Pareto optimal if there is no other solution that is not worse concerning all the objectives considered and is strictly better on at least one of them. The set of Pareto optimal solutions, called the Pareto front, can be a starting point for a multi-objective optimization problem. However, the Pareto front may contain an excessive number of solutions that must be examined one by one by the DM in search of the preferred one. In this perspective, interactive multiobjective optimization methods  provide the DM with a small set of the most preferred solutions seem more suitable for the considered territorial planning problems. Handling multiobjective optimization problems is not straightforward \citep{ehrgott2000survey}. Several methods have been proposed in the literature, such as exact methods \citep{mavrotas2015robustness}, heuristic and metaheuristic methods \citep{ehrgott2008hybrid}. In some cases, the whole Pareto front or part of it can be built \citep{zhou2018multi}, but this makes it difficult for the DM to decide a posteriori which set of facilities is best suitable for him. In any case, determining first the whole Pareto front and after selecting the most desirable solution or looking for the DM's most preferred solutions from the beginning, the preferences of the DM have to be taken into account appropriately. 

In our approach, we propose to integrate a new elicitation algorithm, the above-mentioned DOR procedure, in an interactive multiobjective optimization method looking for the most preferred solution in a procedure alternating discussion and calculation phases as follows.  First, we compute some solutions to an optimization problem. Second, we present  those solutions to the DM and ask him to rank and pairwise compare them in terms of intensity of preferences using the deck of the cards method \citep{figueira2002determining,abastante2020introduction}. Third, the score of the reference solutions obtained through the deck of the cards method is used as  a basis  for  an ordinal  regression procedure \citep{jacquet1982assessing} that,  to  take  into account the interaction among the criteria, is based on a value function expressed in terms of a Choquet Integral \citep{choquet1953theory, Grabisch1996}.  Fourth, the value function so obtained is used to find candidate solutions to the multiobjective optimization problem.  The new feasible solutions can be proposed to the DM, who is again asked to comment on those solutions and rank and pairwise compare them. The process continues  until the DM is satisfied with one of the proposed solutions. The whole iterative process can be supported by using appropriate graphical charts to illustrate the solutions obtained to support the DM throughout the process. As we use a value function that aggregates criteria to evaluate solutions of the multiobjective optimization problem, in the following, we shall use the terms criterion and objective as equivalent.

Our approach presents several advantages:
\begin{itemize}
    \item The DM can participate in the decision making process by expressing his preferences in an easy way thanks to the use of the deck of the cards method.
    \item The deck of the cards method is applied for eliciting the DM's preferences on the solutions of an optimization model instead of being used for expressing more abstract judgments on the importance and interaction of criteria. In this way, the cognitive burden of the DM is reduced, allowing him to  directly comment on some ``feasible" plans and making the process easier and more similar to what happens in reality.
    \item The ordinal regression model permits eliciting the preferences of the DM and dealing with interaction among criteria.
    \item The DM can iteratively build the solutions together with the analyst returning on his preferences at every step of the process.
    \item The whole process is transparent and simple for the DM and supplies arguments to explain the selected solutions to other stakeholders to arrive at a participated decision.
\end{itemize}

We applied the above methodology to support a sustainable territorial decision making process by means of which the following three questions should be answered  in the context of the so-called \textit{space-time model}  proposed by \citep{barbati2020general}:
\begin{itemize}
\item[$\alpha$)]  What facilities need to be selected when planning for a territory? 
\item[$\beta$)]	Where should we locate these facilities?   
\item[$\gamma$)] When those facilities should be activated? 
\end{itemize}
In complex real-world decision problems, all these three questions must be considered simultaneously. Indeed, it is sporadic, particularly in large multi-million euro planning procedures, that a developer can do everything in one shot. Furthermore, administrations and developers  are increasingly pushing for a careful study of the scheduling of the interventions in the plan due to several restrictions or constraints, such as budget constraints, that need to be considered. Several optimization models deal with only some aspects of the planning of territorial projects, answering only one of the three above questions $\alpha)-\beta)-\gamma)$, e.g. \cite{tervonen2017modeling} for $\alpha$),  \cite{farahani2019or} for $\beta$) and \cite{le2020scheduling} for $\gamma$), or a combination of them \cite{sarnataro2021portfolio} for $\beta)$ and $\gamma)$. Instead, adopting the space-time model, we develop  an approach that supports the strategic decision of answering all of the three $\alpha)-\beta)-\gamma)$ questions simultaneously. 

We tested the methodology for setting up an Ecovillage in the  Piedmont Region of Italy. 
The use of the space-time model and of the interactive procedure  is particularly indicated for such a problem for the following reasons: 
\begin{itemize}
\item The DM can realize that the Ecovillage should be treated as a whole system in which the choices related to the facilities to install, to their location and to their timing are related each others in a common overall perspective strategy for which the space-time model seems the most natural methodological scheme.
\item The DM can verify that budget and technical requirements imply a temporal distribution in the activation of the involved facilities.
 \item The DM can recognize that in the setup of an Ecovillage,  a variety of criteria have to be taken into consideration because of its own characteristic of being a self-sufficient village and not a mere profitable investment. Those criteria can also be quite different from more classical criteria in decisions related to conventional touristic structures. In the specific case of an Ecovillage, these criteria present a certain interaction between them that has to be taken adequately into consideration and, in this perspective, we generalize the space-time model to the consideration of the interaction between criteria. Moreover, the weights and the interaction of the considered criteria,  their definition, and their interaction is not always clearly intelligible even for the DMs. In this sense, the use of the DOR methodology permitted an easily understandable indirect preference elicitation procedure because, in this way, the DM was asked to compare some feasible plans comprehensively through the effortless and straightforward procedure of the deck of  the cards method.  Instead, a direct preference elicitation procedure requiring the DM abstract information in terms of importance and interaction of considered criteria would be much more complex and cognitively demanding, with the risk of obtaining not enough reliable results. 
 \item Finally, the introduction of an interactive multiobjective methodology, also thanks  to the newly introduced DOR elicitation procedure,   helps to make a participated decision, taking into account the point of view of the  DM in order to guarantee openness and transparency to the public, in the general perspective of a decision model co-constructed by the analyst with the DM \citep{roy1993decision}.
\end{itemize}
The paper is structured in four parts. After the introduction, Section \ref{sec:proposed methodology} outlines the methodological approach, with a specific focus on the DOR elicitation procedure, while Section \ref{casestudy} describes the real world problem. Section \ref{results} describes the interaction process conducted with the DM and the results obtained, while the last section deals with the conclusions and possible research developments.

\section{The proposed methodology}\label{sec:proposed methodology}

In this Section, we present the methodology we adopted to handle the territorial planning problem. This methodology is based on three main elements: \begin{itemize}
    \item The space time model represents the territorial planning problem in terms of a multiobjective optimization problem answering the three questions of what to do, where to do, and when to do.
   \item  The new DOR elicitation procedure allows the DM's preferences to be collected in an easy and comprehensible manner.
\item An interactive multiobjective optimization procedure permits to search for the most preferred solutions by integrating the above elicitation procedure in  the territorial planning problem expressed in terms of the space time model.  
\end{itemize}
In this perspective, a special role is played by the new elicitation procedure DOR to represent the preferences of DM through a value function. In particular, it is based on a combination of the deck of the cards method \citep{figueira2002determining} in the formulation proposed in \citep{abastante2020introduction} (SRF-II) with an ordinal regression method based on the Choquet integral \citep{marichal2000determination}. We shall use this elicitation procedure to handle an optimization problem formulated in terms of the space-time model \citep{barbati2020general}. However, it has an autonomous interest and it can be used to induce the preference parameters of the Choquet integral as well as of others multicriteria aggregation procedures such as the simplest weighted sum or the piecewise additive utility function considered in UTA method \citep{jacquet1982assessing}. In the following, we formally introduce the space time model and the new elicitation method DOR, while in the last subsection, we summarise the steps of the methodology.

\subsection{The space time model}\label{sec:model}

Let us consider a set of facilities $I=\{1, \ldots, i, \ldots, n\}$. For each facility $i \in I$ we define a set of potential locations $L(i)=\{1(i), \ldots, l(i), \ldots, n(i)\}$. A facility can be assigned to a location at different time epochs $T=\{0,  \ldots, t, \ldots, p\}$. Each facility is evaluated with respect to a set of criteria $G=\{g_j, j \in J\}$ and $J=\{1, \ldots,m\}$. The evaluation of facility $i\in I$ activated in location $l \in L(i)$ with respect to the criterion $g_j\in J$ is denoted by $y_{ijl} \in \mathbb{R}^+$. For the sake of simplicity, without the loss of generality, we suppose that all criteria $g_j\in G$ are of the gain type,  that is, the greater $y_{ijl}$, the better the evaluation of facility $i\in I$ on criterion $g_j \in J$ in location $l \in L(i)$. 

For each time epoch $t \in T$ a discount factor $v(t)$, with $ 0 \le v(t) \le 1$ and $v$ being a non increasing function of $t$, is defined in order to discount the evaluation of performances $y_{ijl}, i \in I, j \in J, l \in L(i)$ in future periods. The values $v(t)$, $t \in T$,  represent the DM's intertemporal preferences. For the sake of simplicity, at first the performances on different criteria are aggregated by abstracting from any consideration on interaction between criteria, so that to make  homogeneous  performances on considered criteria $g_j$ we take into account  weights $w_j \ge 0, j=1, \ldots,m$, which  permit to define an overall value of each plan by summing up the weighted discounted single criterion performances $w_j \cdot y_{ijl} \cdot v(t)$. Each facility $i \in I$ has also a cost $c_{il}\in\mathbb{R}^+$.   The available budget for each period $t \in T$ is denoted by $B_t$.

The following decision variables can be considered to define the adopted plan $\mathbf{x}$:

$$x_{ilt}=
\left\{
\begin{array}{rl}
1, & \hbox{ if facility $i \in I$ is installed in location $l \in L(i) $ in period $t \in T-\{0\}$}; \\
0, & \hbox{ otherwise.}
\end{array} \right.
$$

For example, having a set of facilities $I=\{1,2\}$, the set of locations $L(1)=\{1,2\}$ and $L(2)=\{1,2,3\}$ a set of time epochs $T=\{0,1,2\}$ we have to consider the following vector of decision variables:
$$\mathbf{x}=[x_{110},x_{111},x_{120},x_{121},x_{210},x_{211},x_{220},x_{221},x_{230},x_{231}].$$
If we have
$$x_{110}=x_{111}=x_{120}=0, x_{121}=x_{230}=1,  x_{210}=x_{211}=x_{220}=x_{221}=x_{231}=0,$$
then the adopted plan consists in allocating facility 1 in its second potential location at period 1 and function 2 in its third potential location 1 at period 0.

In case no interaction between criteria is considered, the overall objective function of the space-time optimization model aggregating all the contributions of all criteria in all location and in all time with respect to a plan $\mathbf{x}$ can be formulated as follows: 
\begin{equation}U(\mathbf{x})=\sum_{i \in I}\sum_{j \in J}\sum_{l \in L(i)}\sum_{t \in T-\{0\}}\sum_{\tau=0}^{t-1} v(t)w_j x_{il\tau}y_{ijl}.\end{equation}\label{objFunc1}

Observe that for each criterion $g_j \in G$ and for each plan $\mathbf{x}=[x_{ilt}]$, it is possible to define the overall contribution of criterion $g_j(\mathbf{x})$ as  
\begin{equation}g_j(\mathbf{x})=\sum_{i \in I}\sum_{l \in L(i)}\sum_{t \in T-\{0\}}\sum_{\tau=0}^{t-1} v(t) x_{il\tau}y_{ijl},\end{equation}\label{objFunc2}
so that we can write
\begin{equation}U(\mathbf{x})=\sum_{j \in J}w_jg_j(\mathbf{x}).\end{equation}\label{objFunc3}
%%\begin{equation}\sum_{i \in I}\sum_{j \in J}\sum_{l \in L(i)}\sum_{t \in T-\{0\}}\sum_{\tau=0}^{t-1} v(\tau)w_j x_{il\tau}y_{ijl} + \sum_{sr \in S}\sum_{t \in T-\{0\}} \gamma_t y^{sr}_{jt}.\end{equation}

Observe that not all 0-1 vectors $\mathbf{x}=[x_{ilt}]$ could be feasible. Indeed, a variety of constraints can be defined, according to the particular application at hand, such as: 

 \begin{enumerate}
     \item 
\textbf{Budget constraints} according to which in each period $t \in T$ the expenses cannot be greater than the available budget $B_t$ increased by the possible  unspent budgets from previous periods: 
\begin{equation}\label{Cons:1}
\sum_{i\in I} \sum_{l \in L(i)} c_{il}  x_{ilt} \leq B_t + \sum_{\tau \in T: \tau <  t} B_\tau - \sum_{\tau \in T: \tau <  t}\sum_{i\in I} \sum_{l \in L(i)} c_{il}  x_{il\tau}, \;\; \forall  t\in T,
\end{equation}
that is, in an equivalent formulation,   
\begin{equation}\label{Cons:1_}
\sum_{\tau \in T: \tau \le  t}\sum_{i\in I} \sum_{l \in L(i)} c_{il}  x_{il\tau} \leq B_t + \sum_{\tau \in T: \tau <  t} B_\tau, \;\; \forall  t\in T,
\end{equation}
which can be interpreted by saying that in each period $t$ the total expenses cannot be greater than the sum of all available budgets until $t$; 

\item \textbf{Activation constraints}, i.e. each facility can be activated at most once
\begin{equation}\label{Cons:2}
 \sum_{l(i)\in L(i), t \in T} x_{ilt} \leq 1, \;\; \forall  i\in I.
\end{equation}

\item \textbf{Exclusion constraints}: it may happen that potential locations for different facilities could be the same. In this case it may be impossible to activate both the facilities; we name such a situation as exclusion.
Let us define the set of exclusions $E=\{1, \ldots, e_k,\ldots,K\}$. To each $e_k \in E$ we associate a pair of facilities $\{i,i^\prime\} \in I$ and a pair of respective potential locations $\{l, l'\} \in L(i), L(q)$, respectively. If the facility $i$ is planned in location $l$ then facility $i^\prime$ can not be located in $l'$ at any period $t \in T$. This can be described with the following constraints:
\begin{equation}\label{Cons:3}
\sum_{t \in T} x_{ilt} +\sum_{t \in T}  x_{i'l't} \leq 1,\;\; \forall (\{i,i^\prime\},\{l, l'\})=e_k \in E.\end{equation}

\item \textbf{Scheduling constraints}:
it may happen that some facilities need to be scheduled earlier or later than some other facilities. For instance, if a facility $i$ needs to be scheduled in a earlier period than a facility $i^\prime$, then the following constraints has to be considered:

\begin{equation}\label{Cons:4}
 x_{ilt} \leq \sum_{\tau=0}^{t-1}  x_{i^\prime l \tau}, \;\; \forall t \in T, \forall l \in L.\end{equation}

\end{enumerate}
Another types of constraints are related to consideration of synergy effects between selected facilities in the objective function of the space-time model. More precisely, we consider the case in which the contribution to the different criteria $g_j \in J$ is boosted when some facilities are implemented conjointly in some ``favourable" locations. With this aim, we define a set of synergies $S=\{s_1,\ldots, s_r,\ldots,s_{\overline{r}}\}$, with $s_r=\{i,i^\prime,l, l'\}$, $i,i^\prime \in I, l\in L(i),l^\prime \in L(i^\prime)$.  The synergy $s_r$ is realized when facility $i$ is located in $l$ and facility $i^\prime$ is located in $l^\prime$. In this case, for the period $t$ in which the synergy is realized,  there is an additional contribution $y^r_{jt}=\sigma_{r}\cdot(y_{ijl}+y_{i^\prime j l^\prime})$, with $\sigma_{r} \ge 0$.
To take into account those synergies in our model, we define for each synergy $s_r=\{i,i^\prime,l, l'\} \in S$ and for each $t \in T$, the auxiliary variables  $\gamma^{r}_t$ as:

$$\gamma^{r}_t=
\left\{
\begin{array}{rl}
1, & \hbox {if facilities $i$ and $i^\prime$ results implemented in $l$ and $l^\prime$  at period $t \in T$ or earlier}; \\
0, & \hbox{ otherwise.}
\end{array} \right.
$$
In simple words, $\gamma^{r}_t=1$ if synergy $s_r \in S$ is realized in $t \in T$, and $\gamma^{r}_t=0$ otherwise, which is ensured by the following constraints:

\begin{equation}\label{Cons:5}
\sum_{\tau \in T:\tau\leq t} x_{i l \tau} + \sum_{\tau \in T:\tau\leq t} x_{i^\prime l^\prime \tau} -1 \leq \gamma^{r}_t,\;\; \forall s_r \in S, \forall t  \in T;\end{equation}
\begin{equation}\sum_{\tau \in T:\tau\leq t} x_{i l \tau}\geq\gamma^{r}_t;\forall s_r \in S, \;\; \forall t  \in T;\end{equation}
\begin{equation}\sum_{\tau \in T:\tau\leq t} x_{i^\prime l^\prime \tau}\geq\gamma^{r}_t;\forall s_r \in S, \;\; \forall t  \in T.\end{equation}

Taking into account the contributions of the synergies between facilities, we can reformulate the objective function of the space-time model as follows:

\begin{equation}U(\mathbf{x})=\sum_{i \in I}\sum_{j \in J}\sum_{l \in L(i)}\sum_{t \in T-\{0\}}\sum_{\tau=0}^{t-1} v(t)w_j x_{il\tau}y_{ijl} + \sum_{s_r \in S}\sum_{t \in T-\{0\}} v(t) w_j \gamma^r_t  y^{r}_{jt}.\end{equation}\label{objFunc4}

Observe that the objective function in the formulation (\ref{objFunc4}) can be expressed in terms of the overall contribution of criteria $g_j \in G$ with respect to plan $\mathbf{x}=[x_{ilt}]$ prorperly redefined  as  
\begin{equation}g_j(\mathbf{x})=\sum_{i \in I}\sum_{l \in L(i)}\sum_{t \in T-\{0\}} v(t) (\sum_{\tau=0}^{t-1}  x_{il\tau}y_{ijl}+\sum_{s_r \in S} \gamma^r_t y^{r}_{jt}),\end{equation}\label{objFunc5}
so that we can write
\begin{equation}U(\mathbf{x})=\sum_{j \in J}w_jg_j(\mathbf{x}).\end{equation}\label{objFunc6}

Let us note that the contribution above could be split in relation to one of more elements such as the facility, the period or the criterion. For instance, one could consider  the overall performance in period $t\in T-\{0\}$ of all facilities $i \in I$, all criteria $j \in J$ and all  locations $l \in L$, that is
$y^{T}_{t}(\mathbf{x})=\sum_{i \in I}\sum_{j \in J}\sum_{l \in L}\sum_{\tau=0}^{t-1} w_j x_{il\tau}y_{ijl}$. This could be helpful to understand how the contributions to the criteria of all the facilities activated evolve during the time.

\subsection{Modeling interaction between criteria through the Choquet Integral}\label{Choquet_integral}

%In addition to synergy between facilities we consider also interaction between criteria generalizing the linear aggregation function considered in the objective function of the space time model in equation (\ref{) 

A further enrichment of the objective function of the space - time model is related to the consideration of interaction between criteria, which can be obtained by generalizing the formulation (\ref{objFunc5}) of $U(\mathbf{x})$ in terms of the Choquet integral \citep{choquet1953theory,Grabisch1996}. With this aim, let us introduce the concept of capacity being a function $\mu:2^G\rightarrow [0,1]$ satisfying the following two properties: 
\begin{itemize}
	\item Normalization: $\mu(\emptyset)=0$ and $\mu(G)=1$,
	\item Monotonicity: for all $A \subseteq B \subseteq G, \mu(A) \le \mu(B)$. 
\end{itemize}
For all $A \subseteq G, \mu(A)$ can be interpreted as a value such that, taking into consideration the formulation (\ref{objFunc6}) of $U(\mathbf{x})$, if for the plan $\mathbf{x}$ we have $g_j(\mathbf{x})=k>0$ for all $g_j \in A$ and $g_j(\mathbf{x})=0$ for all $g_j \notin A$, then $U(\mathbf{x})=k\cdot \mu(A)$.  
Given a plan $\mathbf{x}=[x_{ilt}]$ and a capacity $\mu$, the objective function of the space-time model can be expressed in terms of Choquet integral of overall contributions $g_j(\mathbf{x}), g_j \in G,$ as follows  
\begin{equation}U(\mathbf{x})=\sum_{j=1}^m\mu(\{g_h \in G: g_h(\mathbf{x}) \ge g_{(j)}(\mathbf{x})\})\cdot[g_{(j)}(\mathbf{x})-g_{(j-1)}(\mathbf{x})]\end{equation}\label{objFunc7}
with $g_{(1)}(\mathbf{x}), \ldots,g_{(m)}(\mathbf{x})$ a reordering of $g_{1}(\mathbf{x}), \ldots,g_{m}(\mathbf{x})$ such that 
$$g_{(1)}(\mathbf{x}) \le \ldots \le g_{(m)}(\mathbf{x}),$$
and  $g_{(0)}(\mathbf{x})=0$.

Let us remember that a capacity is additive if for all $A,B \subseteq G$ such that $A \cap B =\emptyset, \mu(A \cup B)=\mu(A)+\mu(B)$. In this case, one can set $\mu(\{g_j\})=w_j$ for all $g_j \in G$ and by normalization and monotonicity properties of $\mu$, we get $w_j \ge 0$ for all $g_j \in G$ and $w_1+\ldots+w_m=1$. Moreover, we obtain also that $U(\mathbf{x})=\sum_{j \in J}w_jg_j(\mathbf{x})$, that is, if the capacity $\mu$ is additive then the Choquet integral formulation of the objective function (\ref{objFunc7}) collapse to the linear formulation of the objective function (\ref{objFunc6}). 
  
In case additivity does not hold, the criteria $g_j \in G$ interact between them. For the sake of simplicity, we consider a specific form of interaction that permits obtaining manageable models however permitting them to represent quite general situations. More precisely we consider a 2-additive capacity \citep{grabisch1997k}, that is a capacity $\mu$ such that
there exists $w_j, j=1,\ldots,m$, and $w_{jj^\prime}, \{j,j^\prime\} \subseteq G$ such that, for all $A \subseteq G$
\begin{equation}\mu(A)=\sum_{g_j\in A}w_j+\sum_{\{g_j,g_j^\prime\}\subseteq A}w_{jj^\prime}\end{equation}\label{2additive}
With respect to 2-additive capacities, normalization and monotonicity properties can be reformulated as
\begin{itemize}
	\item Normalization: $\sum_{g_j\in A}w_j+\sum_{\{g_j,g_j^\prime\}\subseteq A}w_{jj^\prime}=1$,
	\item Monotonicity: $w_j\ge 0$ for all $g_j \in G$ and 
\begin{equation}\label{cons:6}
\displaystyle w_j+\sum_{g_{j^\prime}\in T}w_{jj^\prime}\geq 0, \;\mbox{for all} \;g_j\in G \;\mbox{and for all} \;T\subseteq G\setminus\{g_j\}, T\neq\emptyset.\\[1mm]
\end{equation}
\end{itemize}
In case $\mu$ is a 2-additive capacity, the Choquet integral can be expressed as follows 
\begin{equation}\label{ChoquetMobius2additive}
U(\mathbf{x})=\displaystyle\sum_{g_j\in G}w_jg_j(\mathbf{x})+\sum_{\{g_j,g_j^\prime\}\subseteq G}w_{jj^\prime}\min\{g_j(\mathbf{x}),g_{j^\prime}(\mathbf{x})\}
\end{equation}\label{objFunc8}

\subsection{The Deck of the cards based Ordinal Regression}\label{sec:modelregr}

To define the objective function of the space-time model in terms of formulation (\ref{ChoquetMobius2additive}) in Section \ref{Choquet_integral}, it is necessary to fix the value of weights $w_j, g_j \in G$, and $w_{j,j^\prime}, \{g_j,g_{j^\prime}\}\subseteq G$. With this aim we propose DOR, being a new ordinal regression procedure that considers intensity of preferences expressed through the deck of the cards method \citep{figueira2002determining,abastante2020multiple}. The procedure has the following steps
\begin{itemize}
	\item A set of feasible plans $\mathbf{X}=\{\mathbf{x_1},\ldots,\mathbf{x_p}\}$ is presented to the DM.
	\item The DM rank orders the feasible plans from the worst to the best with possible ex-aequo, in $r, r\le p$, equivalence classes $C_1,\ldots,C_r$, so that $C_1$ contains the plans $\mathbf{x}$ considered the worst, $C_r$ contains the plans $\mathbf{x}$ considered the best, and, in general if the plan $\mathbf{x}$ is contained in the equivalence class $C_s$ and the plan $\mathbf{x}^\prime$ is contained in the equivalence class $C_{s^\prime}$ with $s^\prime > s$, then $\mathbf{x}^\prime$ is preferred to $\mathbf{x}$. In particular, a DM is given a set of cards each one representing a plan $\mathbf{x}\in\mathbf{X}$ and the DM is ordering these cards in agreement with the expressed preferences.
		\item The DM puts a certain number of blank cards $e_s, s=1,\ldots,p-1$, between the cards representing the plans $\mathbf{x}$  in the equivalence class $C_s$  and the cards representing the plans  in the equivalence class $C_{s+1}$, in a way that the larger the number of blank cards the greater the difference of preferences between the overall value $\nu(\mathbf{x})$ assigned to the related plans; the DM has also the possibility to put $e_0$ blank cards between a ``zero level'' and the equivalence class $C_1$;
					\item An evaluation $\nu(\mathbf{x})=v_s, s=1,\ldots,p$, is assigned to each plan $\mathbf{x}\in C_s$ applying the following rule
					$$v_s=v_{s-1}+e_{s-1}+1$$
					so that 
					$$v_s=\sum_{z=0}^{s-1} (e_z+1)=\sum_{z=0}^{s-1} e_z +s$$ 
\item The weights $w_j, g_j \in G$, and $w_{jj^\prime}, \{g_j,g_{j^\prime}\}\subseteq G$ of the 2-additive capacity are obtained by minimizing the sum of positive and negative deviations $\sigma^+(\mathbf{x})$ and $\sigma^-(\mathbf{x})$, $\mathbf{x}\in\mathbf{X}$, between the evaluations $U(\mathbf{x})$ assigned by the Choquet integral and the evaluations $\nu(\mathbf{x})$ assigned through the deck of the cards method, properly scaled through a multiplicative positive costant $k$. With this aim, one has to solve the following LP problem  in the variables $w_j, j=1,\ldots,m$, $w_{jj^\prime}, \{g_j,g_j^\prime\}\subseteq G$,$\sigma^+(\mathbf{x})$,$\sigma^-(\mathbf{x})$ and $k$: 
\begin{equation}\label{regressionmodel}
\begin{array}{l}
\min\sum_{\mathbf{x}\in\mathbf{X}}\sigma^+(\mathbf{x})+\sigma^-(\mathbf{x})\\[1mm]
 \mbox{subject to}\\[1mm]
\left.
\begin{array}{l}
U(\mathbf{x})=\displaystyle\sum_{g_j\in G}w_jg_j(\mathbf{x})+\sum_{\{g_j,g_j^\prime\}\subseteq G}w_{jj^\prime}\min\{g_j(\mathbf{x}),g_{j^\prime}(\mathbf{x})\} \;\; \mbox{for all } \mathbf{x}\in\mathbf{X},\\[1mm]
\sum_{g_j\in G}w_j+\sum_{\{g_j,g_j^\prime\}\subseteq G}w_{jj^\prime}=1,\\[1mm]
%\displaystyle\sum_{j=1}^{m}w_j=1,\\[1mm]
w_j\geqslant 0,\;\mbox{for all}\;j=1,\ldots,m, \\[1mm]
w_j+\sum_{g_{j^\prime}\in T}w_{jj^\prime}\geqslant 0, \;\mbox{for all} \;g_j\in G \;\mbox{and for all} \;T\subseteq G\setminus\{g_j\}, T\neq\emptyset,\\[1mm]
U(\mathbf{x})-\sigma^+(\mathbf{x})+\sigma^-(\mathbf{x})=k\cdot \nu(\mathbf{x})\;\; \mbox{for all } \mathbf{x}\in\mathbf{X},\\[1mm]
k\ge 0,\\[1mm]
\sigma^+(\mathbf{x}) \geqslant 0, \sigma^-(\mathbf{x})\geqslant 0 \mbox{ for all } \mathbf{x}\in\mathbf{X}.\\[1mm]
\end{array}
\right\}
\end{array}
\end{equation}
\end{itemize}

Let us observe that the above LP problem can be applied to any form of value function $U(\mathbf{x})$ and not necessarily to a value function formulated in terms of Choquet integral. Indeed, one can see the above LP problem as follows

\begin{equation}\label{regressionmodel_1}
\begin{array}{l}
\min\sum_{\mathbf{x}\in\mathbf{X}}\sigma^+(\mathbf{x})+\sigma^-(\mathbf{x})\\[1mm]
 \mbox{subject to}\\[1mm]
\left.
\begin{array}{l}
E_{value \; function},\\[1mm]
E_{deck \; of \; the \; card \; basis}\\[1mm]
\end{array}
\right\}
\end{array}
\end{equation}

with

\begin{equation}\label{regressionmodel_2}
\begin{array}{l}
%\min\sum_{\mathbf{x}\in\mathbf{X}}\sigma^+(\mathbf{x})+\sigma^-(\mathbf{x})\\[1mm]
 %\mbox{subject to}\\[1mm]
\left.
\begin{array}{l}
U(\mathbf{x})=\displaystyle\sum_{g_j\in G}w_jg_j(\mathbf{x})+\sum_{\{g_j,g_j^\prime\}\subseteq G}w_{jj^\prime}\min\{g_j(\mathbf{x}),g_{j^\prime}(\mathbf{x})\} \\[1mm]
\sum_{g_j\in G}w_j+\sum_{\{g_j,g_j^\prime\}\subseteq G}w_{jj^\prime}=1,\\[1mm]
%\displaystyle\sum_{j=1}^{m}w_j=1,\\[1mm]
w_j\geqslant 0,\;\mbox{for all}\;j=1,\ldots,m, \\[1mm]
w_j+\sum_{g_{j^\prime}\in T}w_{jj^\prime}\geq 0, \;\mbox{for all} \;g_j\in G \;\mbox{and for all} \;T\subseteq G\setminus\{g_j\}, T\neq\emptyset\\[1mm]
\end{array}
\right\}E_{value \; function}
\end{array}
\end{equation} 
and
\begin{equation}\label{regressionmodel_3}
\begin{array}{l}
%\min\sum_{\mathbf{x}\in\mathbf{X}}\sigma^+(\mathbf{x})+\sigma^-(\mathbf{x})\\[1mm]
 %\mbox{subject to}\\[1mm]
\left.
\begin{array}{l}
U(\mathbf{x})-\sigma^+(\mathbf{x})+\sigma^-(\mathbf{x})=k\cdot \nu(\mathbf{x})\;\; \mbox{for all } \mathbf{x}\in\mathbf{X},\\[1mm]
k\ge 0,\\[1mm]
\sigma^+(\mathbf{x}) \geqslant 0, \sigma^-(\mathbf{x})\geqslant 0 \mbox{ for all } \mathbf{x}\in\mathbf{X}.\\[1mm]
\end{array}
\right\}E_{deck \; of \; the \; card \; basis}
\end{array}
\end{equation}

In this perspective, if the considered aggregation function is simply a weighted sum, that is 
$$U(\mathbf{x})=\displaystyle\sum_{g_j\in G}w_jg_j(\mathbf{x})$$
then the set of constraints $E_{value \; function}$ is formulated as follows:
\begin{equation}\label{regressionmodel_4}
\begin{array}{l}
%\min\sum_{\mathbf{x}\in\mathbf{X}}\sigma^+(\mathbf{x})+\sigma^-(\mathbf{x})\\[1mm]
 %\mbox{subject to}\\[1mm]
\left.
\begin{array}{l}
U(\mathbf{x})=\displaystyle\sum_{g_j\in G}w_jg_j(\mathbf{x}),\\[1mm]
\displaystyle\sum_{j=1}^{m}w_j=1,\\[1mm]
w_j\geqslant 0,\;\mbox{for all}\;j=1,\ldots,m. \\[1mm]
\end{array}
\right\}E_{linear \; value \; function}
\end{array}
\end{equation} 

Using the additive piecewise linear value function proposed in the seminal paper \citep{jacquet1982assessing} of the ordinal regression, taking into account criteria $g_j\in G$ assigning each alternative $\mathbf{x}$ a value $g_j(\mathbf{x})$ in the interval $[y_j^0,y_j^{\gamma_j}]$ divided in intervals 
$$[y_j^0,y_j^1],\ldots,[y_j^r,y_j^{r+1}],\ldots,[y_j^{\gamma_j-1},y_j^{\gamma_j}]$$
the set of constraints $E_{value \; function}$ is formulated as follows:
\begin{equation}\label{regressionmodel_5}
\begin{array}{l}
%\min\sum_{\mathbf{x}\in\mathbf{X}}\sigma^+(\mathbf{x})+\sigma^-(\mathbf{x})\\[1mm]
 %\mbox{subject to}\\[1mm]
\left.
\begin{array}{l}
U(\mathbf{x})=\displaystyle\sum_{g_j\in G}u_j(g_j(\mathbf{x}))\\[1mm]
u_j(g_j(\mathbf{x}))=u_j(y^r_j)+\frac{g_j(\mathbf{x})-y^r_j}{y^{r+1}_j-y^r_j}[u_j(y^{r+1}_j)-u_j(y^{r}_j)]\\[1mm]
\;\;\;\;\;\;\;\;\;\;\;\; \mbox{for } g_j(\mathbf{x}) \in [y^{r}_j,y^{r+1}_j],r=0,\ldots, \gamma_j-1,\\[1mm] 
u_j(y^{r+1}_j) \ge u_j(y^r_j) \mbox{ for all} j=1,\ldots,m, r=0,\ldots, \gamma_j-1,\\[1mm]
u_j(y^{0}_j)=0\mbox{ for all } j=1,\ldots,m,\\[1mm]
\displaystyle\sum_{g_j\in G}u_j(y^{\gamma_j})=1.\end{array}
\right\}E_{piecewise\; linear\; value \; function}
\end{array}
\end{equation} 
Observe that the normalization constraint
$$\displaystyle\sum_{g_j\in G}u_j(y^{\gamma_j})=1$$
can be substituted by any constraint 
$$\displaystyle\sum_{g_j\in G}u_j(y^{\gamma_j})=\overline{U}, \overline{U} \in \mathbb{R}^+.$$
For example, in the didactic example of the next section, for the sake of a greater expressivity, we take $\overline{U}=100$.

Some final remarks are useful at the end of this section:
\begin{itemize}
	\item The above proposed DOR elicitation procedure has a general interest not limited to multiobjective optimization problems or decision problems related to territorial planning. In fact, DOR can be used for any MCDA problem taking advantage of the two largely adopted approaches of ordinal regression and the deck of the cards method.
	\item The selection of the analytical form of the value function, that is, taking into account the three cases we considered above, the weighted sum, the Choquet integral, or the additive piecewise linear value function, depends on the specific nature of the decision problem. In general, we can say that: 
		\begin{itemize}
		\item If there is an interest in working with a decision model as simple as possible, the weighted sum should be selected.
		\item If  interactions between criteria have to be taken into account, as is the case for the case study we are taking into consideration in the real-world application, the Choquet integral seems the most adequate formulation of the value function.
		\item If there is an interest in considering how the contribution to the value function of each criterion changes  from one level to the other, the additive piecewise linear value function should be selected. 
	\end{itemize}
	\item The DOR methodology can be applied to any type of aggregation function. The weighted sum, the Choquet integral and the piecewise linear value functions that we considered above are only the most commonly adopted formulations of value functions. For example, in view of the possibility of considering the interaction between criteria without the necessity to express the evaluations of all the criteria on the same scale as is the case of the Choquet integral, it would be interesting to apply the DOR approach to the  enriched additive value function proposed in \citep{greco2014robust}.
\end{itemize}

 %U(\mathbf{x})=\displaystyle\sum_{g_j\in G}w_jg_j(\mathbf{x})+\sum_{\{g_j,g_j^\prime\}\subseteq G}w_{jj^\prime}\min\{g_j(\mathbf{x}),g_{j^\prime}(\mathbf{x})\} \;\; \mbox{for all } \mathbf{x}\in\mathbf{X},\\[1mm]
%U(\mathbf{x})-\sigma^+(\mathbf{x})+\sigma^-(\mathbf{x})=k\cdot \nu(\mathbf{x})\;\; \mbox{for all } \mathbf{x}\in\mathbf{X},\\[1mm]
%\sum_{g_j\in G}w_j+\sum_{\{g_j,g_j^\prime\}\subseteq G}w_{jj^\prime}=1,\\[1mm]
%\displaystyle\sum_{j=1}^{m}w_j=1,\\[1mm]
%w_j\geqslant 0,\;\mbox{for all}\;j=1,\ldots,m, \\[1mm]
%w_j+\sum_{g_{j^\prime}\in T}w_{jj^\prime}\geq 0, \;\mbox{for all} \;g_j\in G \;\mbox{and for all} \;T\subseteq G\setminus\{g_j\}, T\neq\emptyset.\\[1mm]
%k\ge 0,\\[1mm]
%\sigma^+(\mathbf{x})+\sigma^-(\mathbf{x})\mbox{for all } \mathbf{x}\in\mathbf{X}.\\[1mm]
%\end{array}
%\right\}E^{linear}_{DM}
%\end{array}
%$$
\subsection{A didactic example}

In this section, we illustrate  with a simple didactic example the procedure to induce a value function by means of the DOR method. Suppose we have six projects $P_1,P_2,P_3,P_4,P_5$ and $P_6$ evaluated on a $[0-100]$ scale with respect to the three criteria economic aspect, $g_1$, social aspect, $g_2$, and environmental aspect, $g_3$, as shown in the following Table \ref{tab:evaluations}.

\begin{table}[!h]
\begin{center}
\caption{Evaluations of projects with respect to considered criteria }\label{tab:evaluations}
		\begin{tabular}{cccc}
	\hline
   Projects  & Economic aspect: $g_1$ & Social aspect: $g_2$ & Environmental aspect: $g_3$\\
	\hline
${ {\tt P_1}}$ & 80 &	50 &	75  \\  
   ${ {\tt P_2}}$& 60 &	60 &	60  \\
${ {\tt P_3}}$& 60 &	80 &	50  \\
${ {\tt P_4}}$ & 70 &	60 &	70  \\
${ {\tt P_5}}$ & 50 &	70 &	60  \\
${ {\tt P_6}}$& 90 &	50 &	40  \\
		\hline
		\end{tabular}
\end{center}
\end{table}

Using the deck of the cards method and taking into consideration a ``zero project'' $P_0$ as a reference of a null value level, the DM orders the projects from the worst $P_{\{1\}}$ to the best $P_{\{6\}}$, with the number of blank cards $e_s$ between the project $P_{\{s-1\}}$ and the following $P_{\{s\}}, s=1,\ldots,6$, written between brackets $[ \; ]$,  as follows:
$$  {{\tt P_0}} \;\; [30] \;\; {{\tt P_5}} \;\; [3] \;\;{{\tt P_2}} \;\; [1] \;\;{{\tt P_3}} \;\; [6] \;\;{{\tt P_6}} \;\; [1]\;\; {{\tt P_4}} \;\; [4] \;\;{{\tt P_1}}$$

Applying the deck of the cards method we assign the following value to each project
\begin{itemize}
	\item $\nu(P_0=[0,0,0])=0$,
	\item $\nu(P_{\{1\}}=P_5=[50,70,60])=\nu(P_0)+e_1+1=31$,
	\item $\nu(P_{\{2\}}=P_2=[60,60,60])=\nu(P_5)+e_2+1=35$,
	\item $\nu(P_{\{3\}}=P_3=[60,80,50])=\nu(P_2)+e_3+1=37$,
	\item $\nu(P_{\{4\}}=P_6=[90,50,40])=\nu(P_3)+e_4+1=44$,
	\item $\nu(P_{\{5\}}=P_4=[70,60,70])=\nu(P_6)+e_5+1=46$,
	\item $\nu(P_{\{6\}}=P_1=[80,50,75])=\nu(P_4)+e_6+1=51$.
\end{itemize}

%\begin{itemize}
	%\item $\nu(P_0=[0,0,0])=0$,
	%\item $\nu(P_{\{1\}}=P_5=[50,70,60])=\nu(P_0=[0,0,0])+e_1+1=31$,
	%\item $\nu(P_{\{2\}}=P_2=[60,60,60])=\nu(P_5=[50,70,60])+e_2+1=35$,
	%\item $\nu(P_{\{3\}}=P_3=[60,80,50])=\nu(P_2=[60,60,60])+e_3+1=37$,
	%\item $\nu(P_{\{4\}}=P_6=[90,50,40])=\nu(P_3=[60,80,50])+e_4+1=44$,
	%\item $\nu(P_{\{5\}}=P_4=[70,60,70])=\nu(P_6=[90,50,40])+e_5+1=46$,
	%\item $\nu(P_{\{6\}}=P_1=[80,50,75])=\nu(P_6=[70,60,70])+e_6+1=51$.
%\end{itemize}

Taking into account a value function $U$ expressed in terms of a weighted sum, the ordinal regression methodology proposed in Section \ref{sec:modelregr} can then be applied to solve the following LP problem in the variables $w_1, w_2, w_3, \sigma^+(\mathbf{P_i}),\sigma^-(\mathbf{P_i}), i=1,\ldots,6$, and $k$: 
\begin{equation}\label{regressionmodel_WS}
\begin{array}{l}
\min\sum_{i=1}^6\sigma^+(\mathbf{P}_i)+\sigma^-(\mathbf{P}_i)\\[1mm]
 \mbox{subject to}\\[1mm]
\left.
\begin{array}{l}
U(\mathbf{P}_i)=w_1g_1(\mathbf{P}_i)+w_2g_2(\mathbf{P}_i)+w_3g_3(\mathbf{P}_i) \;\; i=1,\ldots,6\\[1mm]
U(\mathbf{P}_i)-\sigma^+(\mathbf{P}_i)+\sigma^-(\mathbf{P}_i)=k\cdot \nu(\mathbf{P}_i),\;\; i=1,\ldots,6,\\[1mm]
w_1+w_2+w_3=1,\\[1mm]
%\displaystyle\sum_{j=1}^{m}w_j=1,\\[1mm]
w_1\geqslant 0,w_2\geqslant 0,w_3\geqslant 0, \\[1mm]
k\ge 0,\\[1mm]
\sigma^+(\mathbf{P}_i) \geqslant 0, \sigma^-(\mathbf{P}_i)\geqslant 0, \;\; i=1,\ldots,6.\\[1mm]
\end{array}
\right\}
\end{array}
\end{equation}
The solution of the LP problem (\ref{regressionmodel_WS}) gives the results shown in Table \ref{Tab:regressionmodel_WS} with the scaling costant $k=1.282$ and the following weights for the considered criteria: $w_1= 0.517, w_2=0.079, w_3=0.404$. The total sum of errors $\sum_{i=1}^6\sigma^+(\mathbf{P}_i)+\sigma^-(\mathbf{P}_i)$ is 8.09.

\begin{table}[!h]
\begin{center}
\caption{Scores assigned to projects by the value function $U$ obtained solving the LP problem  (\ref{regressionmodel_WS})}\label{Tab:regressionmodel_WS}
		\begin{tabular}{cccccc}
	\hline
   Projects  & $U(\mathbf{P}_i)$ & $\nu(\mathbf{P}_i)$ & $k\cdot\nu(\mathbf{P}_i)$& $\sigma^+(\mathbf{P}_i)$ & $\sigma^-(\mathbf{P}_i)$ \\
	\hline
${ {\tt P_1}}$ & 75.62 &	51 &	75.62 & 0&0\\  
   ${ {\tt P_2}}$& 60 &	35 &	55.12  &4.8&0\\
${ {\tt P_3}}$& 57.53 &	37 &	57.68  &0&0.15\\
${ {\tt P_4}}$ & 69.21 &	46 &	69.21  &0&0\\
${ {\tt P_5}}$ & 55.61 &	31 &	52.55  &3.06&0\\
${ {\tt P_6}}$& 66.65 &	44 &	66.65  &0&0\\
		\hline
		\end{tabular}
\end{center}
\end{table}

If one wants to take into account a value function expressed in terms of piecewise linear function dividing the interval $[0,100]$ of possible values assigned by criteria $g_1,g_2,g_3$ in the intervals 
$$[0,50],[50,75],[75,100],$$
the following LP problem in the variables $u_j(0), u_j(50), u_j(75), u_j(100), j=1,2,3, \sigma^+(\mathbf{P_i}),\sigma^-(\mathbf{P_i}), i=1,\ldots,6$, and $k$ has to be solved: 

\begin{equation}\label{regressionmodel_PLU}
\begin{array}{l}
\min\sum_{i=1}^6\sigma^+(\mathbf{P}_i)+\sigma^-(\mathbf{P}_i)\\[1mm]
 \mbox{subject to}\\[1mm]
\left.
\begin{array}{l}
U(\mathbf{P}_i)-\sigma^+(\mathbf{P}_i)+\sigma^-(\mathbf{P}_i)=k\cdot \nu(\mathbf{P}_i),\;\; i=1,\ldots,6,\\[1mm]
U(\mathbf{P}_i)=\displaystyle\sum_{g_j\in G}u_j(g_j(\mathbf{P}_i)),\\[1mm]
u_j(g_j(\mathbf{P}_i))=u_j(y^r_j)+\frac{g_j(\mathbf{P}_i)-y^r_j}{y^{r+1}_j-y^r_j}[u_j(y^{r+1}_j)-u_j(y^{r}_j)]\;\mbox{for } g_j(\mathbf{P}_i) \in [y^{r}_j,y^{r+1}_j]\\[1mm] 
u_j(50) \ge u_j(0),  j=1,2,3,\\[1mm]
u_j(75) \ge u_j(50),  j=1,2,3,\\[1mm]
u_j(100) \ge u_j(75),  j=1,2,3,\\[1mm]
u_j(0)=0,  j=1,2,3,\\[1mm]
u_1(100)+u_2(100)+u_3(100)=100,\\[1mm]
k\ge 0,\\[1mm]
\sigma^+(\mathbf{P}_i) \geqslant 0, \sigma^-(\mathbf{P}_i)\geqslant 0 \;\; i=1,\ldots,6.\\[1mm]
\end{array}
\right\}
\end{array}
\end{equation}
					
The solution of the LP problem (\ref{regressionmodel_PLU}) gives the marginal value function determined by the values $u_j(0), u_j(50), u_j(75), u_j(100), j=1,2,3$,  shown in Table \ref{regressionmodel_PLU_basic} with the scaling costant $k=1.11$. The projects $\textbf{P}_i,i=1,\ldots,6$, receive the evaluations shown in Table \ref{regressionmodel_PLU_basic}\label{tab:PLU results}.  The total sum of errors $\sum_{i=1}^6\sigma^+(\mathbf{P}_i)+\sigma^-(\mathbf{P}_i)$ is equal to 0. Observe that in the LP problem (26), through the constraint $$u_1(100)+u_1(100)+u_1(100)=100$$
we fixed $\overline{U}=100$.

\begin{table}[!h]
\begin{center}
\caption{Reference values defining the piecewise additive value function $U$ obtained solving the LP problem  (\ref{regressionmodel_PLU})}\label{regressionmodel_PLU_basic}
		\begin{tabular}{ccccc}
	\hline
     & $u_j(0)$ & $u_j(50)$ & $u_j(75)$& $u_j(100)$ \\
	\hline
Economic & 0 &	31.48 &	47.22 & 64.81\\  
   Social & 0 &	0 &	0.10  &0.20 \\
Environmental & 0 &	0 &	14.81  &14.81\\
		\hline
		\end{tabular}
\end{center}
\end{table}

\begin{table}[!h]
\begin{center}
\caption{Scores assigned to plans by the value function $U$ obtained solving the LP problem  (\ref{regressionmodel_PLU})}\label{tab:PLU results}
		\begin{tabular}{cccccc}
	\hline
   Projects  & $U(\mathbf{P}_i)$ & $\nu(\mathbf{P}_i)$ & $k\cdot\nu(\mathbf{P}_i)$& $\sigma^+(\mathbf{P}_i)$ & $\sigma^-(\mathbf{P}_i)$ \\
	\hline
${ {\tt P_1}}$ & 65.56 &	51 &	65.56 & 0&0\\  
   ${ {\tt P_2}}$& 47.78 &	35 &	47.78  &0 &0\\
${ {\tt P_3}}$& 50.00 &	37 &	50.00  &0&0\\
${ {\tt P_4}}$ & 60.00 &	46 &	60.00  &0&0\\
${ {\tt P_5}}$ & 45.56 &	31 &	45.56  &0&0\\
${ {\tt P_6}}$& 57.78 &	44 &	57.78  &0&0\\
		\hline
		\end{tabular}
\end{center}
\end{table}

Finally, taking into account a value function expressed in terms of Choquet integral the following LP problem \ref{regressionmodel_Choquet} has to be solved in the variables $w_1, w_2, w_3, w_{12},
 w_{23}, w_{13},\sigma^+(\mathbf{P_i}),\sigma^-(\mathbf{P_i}), i=1,\ldots,6$, and $k$:

\begin{equation}\label{regressionmodel_Choquet}
\begin{array}{l}
\min\sum_{i=1}^6\sigma^+(\mathbf{P}_i)+\sigma^-(\mathbf{P}_i)\\[1mm]
 \mbox{subject to}\\[1mm]
\left.
\begin{array}{l}
U(\mathbf{P}_i)-\sigma^+(\mathbf{P}_i)+\sigma^-(\mathbf{P}_i)=k\cdot \nu(\mathbf{P}_i)\;\; i=1,\ldots,6,\\[1mm]
U(\mathbf{P}_i)=w_1g_1(\mathbf{P}_i)+w_2g_2(\mathbf{P}_i)+w_3g_3(\mathbf{P}_i)+\\[1mm]
+w_{12}min\{g_1(\mathbf{P}_i),g_2(\mathbf{P}_i)\}+w_{13}min\{g_1(\mathbf{P}_i),g_3(\mathbf{P}_i)\}+w_{23}min\{g_2(\mathbf{P}_i),g_3(\mathbf{P}_i)\},\\[1mm]
w_1+w_2+w_3+w_{12}+w_{23}+w_{13}=1,\\[1mm]
w_j+\sum_{g_{j^\prime}\in T}w_{jj^\prime}\geq 0, \;\mbox{for all} \;g_j\in \{g_1,g_2,g_3\} \;\mbox{and for all} \;T\subseteq \{g_1,g_2,g_3\}\setminus\{g_j\}, T\neq\emptyset,\\[1mm]
k\ge 0,\\[1mm]
\sigma^+(\mathbf{P}_i) \geqslant 0, \sigma^-(\mathbf{P}_i)\geqslant 0, \;\; i=1,\ldots,6.\\[1mm]
\end{array}
\right\}
\end{array}
\end{equation}

The solution of the LP problem (\ref{regressionmodel_Choquet}) gives $w_1=0.52, w_2=0.08, w_3=0.10, w_{12}=0,
 w_{23}=0.31, w_{13}=0$ with the project $\textbf{P}_i,i=1,\ldots,6$, receiving the evaluations shown in Table \ref{tab:regressionmodel_Choquet} and the total sum of errors $\sum_{i=1}^6\sigma^+(\mathbf{P}_i)+\sigma^-(\mathbf{P}_i)$ being equal to 5.03.

\begin{table}[!h]
\begin{center}
\caption{Scores assigned to projects by the value function $U$ obtained solving the LP problem  (\ref{regressionmodel_Choquet})}\label{tab:regressionmodel_Choquet}
		\begin{tabular}{cccccc}
	\hline
   Projects  & $U(\mathbf{P}_i)$ & $\nu(\mathbf{P}_i)$ & $k\cdot\nu(\mathbf{P}_i)$& $\sigma^+(\mathbf{P}_i)$ & $\sigma^-(\mathbf{P}_i)$ \\
	\hline
${ {\tt P_1}}$ & 75.62 &	51 &	75.62 & 0&0\\  
   ${ {\tt P_2}}$& 60.00 &	35 &	55.12  &4.88 &0\\
${ {\tt P_3}}$& 57.53 &	37 &	57.68  &0&0.15\\
${ {\tt P_4}}$ & 69.21 &	46 &	69.21  &0&0\\
${ {\tt P_5}}$ & 52.55 &	31 &	52.55  &0&0\\
${ {\tt P_6}}$& 66.65 &	44 &	66.65  &0&0\\
		\hline
		\end{tabular}
\end{center}
\end{table}

\subsection{Summary of the steps}

In the following, we propose a summary of the steps for the proposed methodology:

\begin{enumerate}
    \item \textbf{Structuring the problem:} The analyst and the DM define the main elements of the problems in terms of objectives/criteria to take into consideration, the facilities, their location, and their evaluations. They also specify the planning horizon  and other characteristics that the plans should have.
    \item \textbf{Identification of potential plans:} The analyst selects some plans to submit to the attention of the DM. This step can be conducted with the definition of some plans obtained, for example, including relevant constraints related to desired characteristics of the plan in the space-time model of Subsection \ref{sec:model} and optimising single criteria.
    \item \textbf{Ranking of the proposed plans and elicitation of the DM preferences:}
    The DM ranks  the proposed plans and compares them with the deck of the card method obtaining an evaluation $\nu(\mathbf{x})$ for each plan $\mathbf{x}$. With the applications of the regression model of Subsection \ref{sec:modelregr} a set of weights $w_j$ for each criterion $g_j$ and a set of interaction coefficients $w_{jj^\prime}, \{g_j,g_{j^\prime}\}\subseteq G$, is derived and a new value function for the space-time model is defined. The DM also comments on the plans obtained and his indications can be introduced as constraints in the subsequent space-time model. 
    \item \textbf{Definition of a new set of plans:}
    Thanks to the application of the  space-time model of Subsection \ref{sec:model} and the value function obtained in the previous step, new plans are generated.  If the DM is satisfied with one of the proposed plans, the procedure stops. Otherwise, we go back to step 3, asking the DM to express his preferences on the newly generated plans, and the procedure is iterated until the DM is satisfied with one of the proposed plans.                                                                                                                                                       
\end{enumerate}

\section{The real-world application}\label{casestudy}
\begin{figure}[htp!]
 \centering
\includegraphics{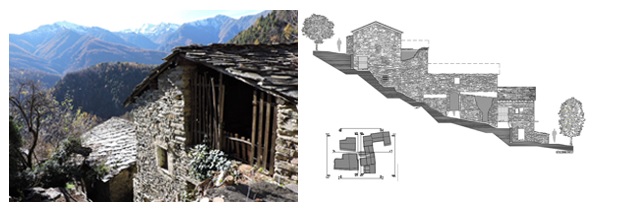}
 \caption{One of the buildings of the ``House of the Sun” and a transformation hypothesis 
(source: libertarea.org)
}\label{fig:Photo}

\end{figure}

The real-world application consists of the implementation of an Ecovillage in Italy. Ecovillages may be considered as rural enterprises that unify sustainable environment-friendly technologies, organic agriculture, and other farming activities and tourism services. 
The Ecovillages represent a type of lifestyle. Based on this philosophy, the Ecovillages are usually designed and built up within the framework of four focuses: the ecologic, social, cultural, and spiritual concepts.
The case under analysis is the hypothesis of the revitalisation of a rural settlement built at the end of the 18th century in dry stone, at an altitude of 1000 meters, abandoned in the 1950s, and located in the mountains about an hour from Turin (the capital of the region). It comprises two small boroughs, the Upper Borough and the Lower Borough, with 11.4 surrounding hectares of woodland (see Figure \ref{fig:Photo}). 
After years of searching and negotiating, a cooperative bought this rural settlement to create an Ecovillage called ``The House of the Sun". Their motto is: ``another world is possible, we are building it... here!". The idea is to be able to restore the relationship with nature and with the environment more harmoniously, through food, furnishings, clothing, and also a whole series of practices, in addition to those working which may be organic farming, even a little more unusual and holistic as the martial arts, yoga, or meditation, rather than shiatsu treatment or tai chi chuan, but also more simply the traditional folk dances to recover the Occitan tradition of these cross-border valleys. This is all part of a dynamic exchange with the territory to reactivate the economic fabric of the valley: the experience of artisans who know how to (still) build with stone and wood and involve those who want to help the cooperative in the enterprise of revitalizing the valley.

Defining the facilities, their locations and the timing of an Ecovillage  is undoubtedly challenging because it is a very particular case of territorial transformation with non-ordinary logic, where, for example, money has a very different value compared to urban transformation contexts where the goal of the developer is to maximize income. There are several peculiar aspects of an Ecovillage that must be considered:
\begin{itemize}
    \item 
	The informal economy plays a fundamental role as one has to consider also exchanges that take place by the social network, without the exchange of money (e.g. barter). This is an important aspect to consider in the location of functions, which for residents will follow non-commercial logic.
\item	There is no certain right or wrong concept while building up an Ecovillage.  What is generally recognized is that a careful and specific design is important for healthy development in the long run. In this sense, Ecovillages use technologies, e.g. passive solar energy design, natural isolation materials, biomass gas converters, etc. 
\item	The social aspect is fundamental. In each Ecovillage, a conscious effort is made to develop the community environment and create a sense of belonging.
\item	The Ecovillage involves the presence of three types of users: i) residents, i.e. people living all year round there; ii) temporary residents, who work in the village for a period ranging from 2 weeks to 6 months, and which have a specific name, WWOOFER (World Wide Opportunities on Organic Farms); iii) guests (as in hotels), who practice a type of tourism with a very strong environmental connotation (eco-tourism).
\end{itemize}
The last point implies that the allocation of services considers which functions could, temporarily or permanently, be used by different types of beneficiaries. For instance, it is possible to imagine that the first two types of users could have similar residential spaces and temporarily share the residents' common areas.  In general, all the spaces have to be created in order to stimulate interactions, protect privacy, and give the opportunity to develop a sense of community.
The decision to use this case study arose from the opportunity to interact with the President of the cooperative owning ``The House of the Sun” (hereinafter defined as the Decision Maker, DM). The strong conviction to create an alternative way of living and working clashed with severe budget constraints. In this sense, applying our new elicitation procedure for handling multiobjective combinatorial optimization problems seemed to be perfect in relation to this real-world problem.

\section{Results and Implementation of the methodology}\label{results}
\subsection{Structuring the problem}

In collaboration with the DM, we structured the problem taking into account the following elements:
\begin{itemize}
    \item 
The set of facilities $I=\{1,\ldots,10\}$, distinguished by those for the residents and those for the tourists (including the WWOOFERs). The facilities to be included concern these two types of users, although the level of interaction between the two can be very strong, especially in the first years of the Ecovillage. Both residents and tourists will need a kitchen, dining room, and rooms, then there is the tailoring/laundry room, the woodworking, and the recreational room (destined for yoga, meditation, martial arts, and dance). More in detail, the facilities with their respective symbols and labels are reported in Table \ref{tab:listprojects}.

\begin{table}[h!]
\centering
 \begin{tabular}{lcc} \hline\hline
 Function & Label&Symbol\\\hline
\textit{Residence for the WWOOFER}&({ {\tt RES-WWO}})&$\includegraphics[height=\myMheight]{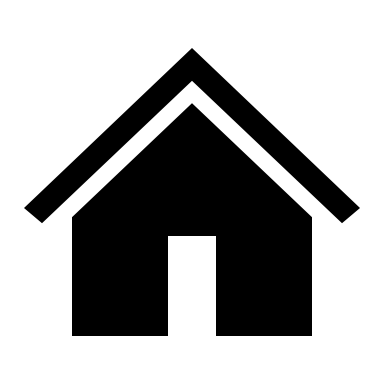}_R$\\
\textit{Kitchen for the WWOOFER}&({ {\tt KIT-WWO}})&$\includegraphics[height=\myMheight]{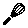}_R$\\
\textit{Refectory for the WWOOFER}&({ {\tt REF-WWO}})&$\includegraphics[height=\myMheight]{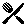}_R$\\
\textit{Guest Rooms}&({ {\tt ROM-GUE}})&$\includegraphics[height=\myMheight]{House_.png}_G$\\
\textit{Guest  Kitchen}&({ {\tt KIT-GUE}})&$\includegraphics[height=\myMheight]{Kitchen_.png}_G$\\
\textit{Guest Dining room}&({ {\tt DIN-GUE}})&$\includegraphics[height=\myMheight]{Dining_.png}_G$\\
\textit{Laboratory 1: tailoring}&({ {\tt TAI-LAB}})&\includegraphics[height=\myMheight]{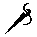}\\
\textit{Laboratory 2: woodworking}&({ {\tt WOO-LAB}})&\includegraphics[height=\myMheight]{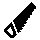}\\
\textit{Recreational room (yoga / meditation martial arts dance)}&({ {\tt ROM-REC}})&\includegraphics[height=\myMheight]{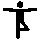}\\
\textit{Main technical room}&({ {\tt ROM-TEC}})&\includegraphics[height=\myMheight]{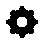}\\

\end{tabular}
 \caption{List of facilities}\label{tab:listprojects}
\end{table}

\item	The sets of locations $L(i)=\{l_1(i),l_2(i)\}$ defining for each facility $i \in I$  the potential location for each facility in the Upper or Lower boroughs (see Figure 2). The two locations are a short distance apart, the upper location is a little larger, but both are in a serious state of disrepair and require extensive renovation. According to the technical and positional characteristics of the different rooms in the buildings in the Upper Borough and the Lower Borough, the  functions can be located only in specific spaces (primarily according to the surfaces needed). All locations are the result of a significant renovation of existing buildings with only possible new construction of a pavilion for recreational activities has been considered.
 In Table \ref{tab:locations} the different spaces are identified with a letter (corresponding to the building) and a number (to distinguish the different rooms located at the different levels of the buildings). 
 \item The  cost $c_{il}$ associated to each location $l \in L(i)$ and to each facility $i \in I$  (see Table \ref{tab:locations}). 
The cost represents an estimation of the implementation costs. In addition to the construction costs indicated in the table, the following items of expenditure have been estimated, and appropriately distributed over the four years considered: Design costs; General expenses; Primary and secondary urbanization charges; Initial costs (purchase of furniture and machinery); Annual running costs. 

\item	The set of periods $T=\{t_0,t_1,t_2,t_3\}$, with $t_0=0,t_1=1,t_2=2,t_3=3$, i.e. we are investigating the possibility that the planning period will last for three years.

\item	The set of criteria $G=\{g_1, g_2, g_3, g_4$\} that have been derived by the analysis of the aims for installing  Ecovillages and they have been extensively discussed with the DM.
More in detail:
\begin{itemize}
    
\item	Environmental ($g_1$): it is the “mother principle” that determines everything else, it is considered the fundamental value that motivates this peculiar choice of life;
\item	Social ($g_2$): it is related to the will to repopulate inland territories (an objective recognized as particularly important at European level and, paradoxically, less at the Italian level), while encouraging urban congestion;
\item	Economic ($g_3$): it considers two main aspects. On  one hand, a principle of self-sustainability with a low environmental impact is a fundamental and structural objective to be pursued; on the other hand, the issue of running a profitable activity related to the eco-tourism;
\item	Cultural ($g_4$): it takes into account how  activities in the area are intertwined with social and cultural themes (e.g. guided socio-hiking, rediscovery of local history, aggregation of schooling, etc.).
\end{itemize}
Theoretically, these four criteria must always be optimised together because the ecological cultural  holistic basic assumption implies the consideration of strong interaction between these four criteria. In this perspective, in consideration of its capacity to model the interaction between criteria,  the Choquet integral model seems the most adequate formulation of the value function $U$ for the decision problem at hand. In Table \ref{tab:criteria}, we can see for each facility $ i \in I$ and for each location $l \in L$ the evaluations $y_{ijl}$ for each criterion $g_j \in G$; those estimates were provided by the expert and agreed with the DM and for the sake of the simplicity are expressed with values between 0 and 100. 

%\begin{figure}[htp!]\label{fig:Mappa}
 %\centering
%\includegraphics{Mappa.jpg}
 %\caption{Some of the possible functions locations in the House of the Sun (source: Libertarea)
%}
 %\label{fig:Mappa}
%\end{figure}

\begin{table}[htb!]
\centering
 \begin{tabular}{cp{3.8cm}cp{3.8cm}c} \hline\hline
 Facilities  & Location $1$ & $c_1$ & Location $2$ & $c_2$   \\
 \hline
   ({ {\tt RES-WWO}}) &{B1, B7, B8, B9, B10, A3, A4, A5, A6, A7}&	212,175 €&	H1, H2, H3, H4, I1, I2, I3, I4, L1, L2 &	233,390 €  \\  
   ({ {\tt KIT-WWO}})&	B4&	26,560 €	&M1&	29,215 €\\
({ {\tt REF-WWO}})&	B3&	15,955 €&	M2&	17,550 €\\
({ {\tt ROM-GUE}}) & {F4, F6, A7, D4, D6, C4, C5, C6}&	185,515 €& {	B1, B7, B8, B9, B10, A3, A4, A5, A6, A7}&	212,175 € \\
({ {\tt KIT-GUE}}) &	C2&	18,235 €&	D3&	30,090 €\\
({ {\tt DIN-GUE}})&	C1&	31,910 €&	D1, E6,E5&	73,800 €\\
({ {\tt TAI-LAB}})&	B6& 14,865 €&	C2&	35,100 €\\
({ {\tt WOO-LAB}}) &	C7&	31,910 €&	F6&	8,720 €\\
({ {\tt ROM-REC}})&	C8&	21,405 €&	Pavillon&	23,545 €\\
({ {\tt ROM-TEC}})&	F5&	13,975 €&	H5&	20,060 €

  \\ \hline\hline
 \end{tabular}
 \caption{Locations of each facility and the associated costs}\label{tab:locations}
\end{table}

\begin{table}[htb!]
\centering
 \begin{tabular}{p{2cm}p{1cm}p{1cm}p{1cm}p{1cm}p{1cm}p{1cm}p{1cm}p{1cm}} \hline\hline
\multirow{2}{*}{Facilities}  & \multicolumn{2}{c|}{Environmental}& \multicolumn{2}{c|}{Social }&\multicolumn{2}{c|}{Economic}&\multicolumn{2}{c}{ Cultural}   \\
 \cline{2-9}
 &$l_1$&$l_2$&$l_1$&$l_2$&$l_1$&$l_2$&$l_1$&$l_2$\\
 \hline
   ({ {\tt RES-WWO}}) &		80	&	80	&	82	&	70	&	40	&	35	&	80	&	80	\\
 ({ {\tt KIT-WWO}})&		80	&	80	&	82	&	70	&	40	&	35	&	80	&	80	\\
({ {\tt REF-WWO}})&		80	&	80	&	82	&	70	&	40	&	35	&	80	&	80	\\
({ {\tt ROM-GUE}})&		60	&	60	&	70	&	0	&	72	&	80	&	70	&	70	\\
({ {\tt KIT-GUE}}) &		55	&	60	&	70	&	70	&	72	&	80	&	65	&	70	\\
({ {\tt DIN-GUE}})&		55	&	60	&	62	&	70	&	72	&	80	&	65	&	70	\\
({ {\tt TAI-LAB}})&		70	&	62	&	43	&	38	&	50	&	50	&	70	&	72	\\
({ {\tt WOO-LAB}}) &		70	&	65	&	45	&	40	&	55	&	65	&	70	&	72	\\
({ {\tt ROM-REC}})&		72	&	60	&	55	&	42	&	55	&	70	&	62	&	78	\\
({ {\tt ROM-TEC}})&		75	&	75	&	35	&	35	&	42	&	48	&	72	&	72	\\

  \\ \hline\hline
 \end{tabular}
 \caption{Criteria evaluations for each facility and for each location}\label{tab:criteria}
\end{table}

\item In terms of characteristics that the plans must have, the DM and the analysts agreed that:
\begin{itemize}
    \item  A facility could be activated only once and only in one location.
   \item  The pairs of facilities ({{\tt RES-WWO}}) and ({{\tt ROM-GUE}}), ({{\tt KIT-GUE}}) and ({{\tt WOO-LAB}}) and ({{\tt ROM-TEC}}) should not be opened in the same location.
    \item The facilities ({{\tt KIT-GUE}}) and ({{\tt DIN-GUE}}) if opened at the same time would cause an increase in the evaluation of the facilities with respect to the considered criteria of $\sigma_{r} =20\%$. 
\end{itemize}
  \end{itemize}  
  Each one of the above requirements was considered in the definition of the plan by means of specific constraints included in the formulation of the space-time model. The plans proposed to the DM were obtained by maximizing some specific value function $U$ as detailed in the following.
  
\subsection{Identification of potential plans}
% \section {Results and Discussion} \label{results}

In order to propose some plans to the DM, the analysts adopted the space-time model introduced in Section \ref{sec:model}. More in detail, we simulated two different scenarios according to two different budgets configurations: \begin{itemize}
    \item $100,000$ Euro in every period $t \in T$, called Budget configuration $B_1$;
    \item $50,000$ Euro in every period $t \in T$, called Budget configuration $B_2$.
\end{itemize} 

In this initial stage, we aggregated the evaluations on the considered criteria by means of a value function $U$ expressed in terms of a weighted sum considering four different weight vectors $\mathbf{w}=[w_1,w_2,w_3,w_4]$ collecting weights $w_j$ for criteria $g_j, j=1,2,3,4$, reported in Table \ref{tab:weightsfirst}. Those initial weights were chosen to represent equal weights or give much more importance to one of the criteria than the others. For this initial stage, we did not consider potential interaction among the criteria, and, consequently, we did not adopt a more complex and sophisticated Choquet integral model, because we just wanted to propose some initial plans to the DM in order to start the discussion.  In other words, we can say that in this first step, we fixed  the interaction coefficients $w_{j,j^\prime}, \{g_j,g_{j^\prime}\}\subseteq G$ equal to zero. 

\begin{table}[htb!]
\centering
 \begin{tabular}{c|cccc}
& $w_1$&$w_2$&$w_3$&$w_4$\\
 \hline\hline
$\mathbf{w^1}$&0.25&0.25&0.25&0.25\\
$\mathbf{w^2}$&0.997&0.001&0.001&0.001\\
$\mathbf{w^3}$&0.001&0.997&0.001&0.001\\
$\mathbf{w^4}$&0.001&0.001&0.997&0.001\\
$\mathbf{w^5}$&0.001&0.001&0.001&0.997\\
 \end{tabular}
 \caption{Selected set of weights for the initial stage }\label{tab:weightsfirst}
\end{table}
To the formulation of the space-time model we added the activation constraints (\ref{Cons:2}) for each facility $i \in I$, and 
the exclusion constraints (\ref{Cons:3}) among the pairs of facilites ({{\tt RES-WWO}}) and ({{\tt ROM-GUE}}), ({{\tt KIT-GUE}}) and ({{\tt WOO-LAB}}) and ({{\tt ROM-TEC}}) according to the DM's preferences. We also defined a discount factor $v(t)=1.10^{-t}$.
 In addition to that, we ran all the scenarios defined above  with the synergy constraint between the facilities ({{\tt KIT-GUE}}) and ({{\tt DIN-GUE}}). If those facilities were opened simultaneously, they would create an additional contribution of 20\% for the four criteria considered. Indeed, during our initial discussion with the DM, he expressed that this synergy would be important, but he would happily discuss plans without a synergy. Therefore, to attain a set of initial plans as different as possible, we simulated all the scenarios with this synergy constraint, identified with $SG_1$, and without the synergy constraint, identified with $SG_2$. In this way, maximizing the value function $U(\mathbf{x})=\sum_{g_j \in G}w_jg_j(\mathbf{x})$ in the different scenarios $(B_r,\mathbf{w^s}, S_k)$ obtained by the combination of the budget $B_r, r=1,2$,  the weight vectors $\mathbf{w^s}, s=1,\ldots,5,$ and presence of synergy constraint $SG=\{SG_1, SG_2\}$ we obtained  20 initial plans. Some of these obtained plans were the same. In addition, to reduce the cognitive burden for the DM, we decided to select only the most representative ones and the ones that presented more differences. In the end, eight different plans $\mathbf{x}_1, \ldots, \mathbf{x}_8$ were presented to the DM as reported in Table \ref{tab:firststrategies} with the first four plans obtained with the budget configuration $B_1$ and the other four plans obtained with budget configuration $B_2$; the symbol  $\times$  means that a particular facility has not been selected, otherwise, the location $l \in L$ and the period $t \in T$ in which the facility is implemented are presented.
 More precisely, the plans selected could be obtained in this way:
 \begin{itemize}
     \item $\mathbf{x}_1$, for budget $B_1$, weights $\mathbf{w^1}$, presence of synergy $SG_1$;
      \item $\mathbf{x}_2$, for budget $B_1$, weights $\mathbf{w^5}$, absence of synergy $SG_2$;
      \item $\mathbf{x}_3$, for budget $B_1$, weights $\mathbf{w^4}$, absence of synergy $SG_2$;
      \item $\mathbf{x}_4$, for budget $B_1$, weights $\mathbf{w^3}$, absence of synergy $SG_2$;
      \item $\mathbf{x}_5$, for budget $B_2$, weights $\mathbf{w^3}$, presence of synergy $SG_1$;
      \item $\mathbf{x}_6$, for budget $B_2$, weights $\mathbf{w^1}$, presence of synergy $SG_1$;
      \item $\mathbf{x}_7$, for budget $B_2$, weights $\mathbf{w^5}$, absence of synergy $SG_2$;
      \item $\mathbf{x}_8$, for budget $B_2$, weights $\mathbf{w^4}$, absence of synergy $SG_2$;

 \end{itemize}
 Note that each plan can be obtained more times for different combinations of parameters as plan $\mathbf{x}_6$ that is the same also in the absence of synergy $G_1$.

\begin{table}[htb!]
\centering
 \begin{tabular}{c|cccccccccc}

  & \scriptsize({ {\tt RES-WWO}}) & \scriptsize({{\tt KIT-WWO}})&\scriptsize({{\tt REF-WWO}})&\scriptsize({{\tt ROM-GUE}}) &
\scriptsize({{\tt KIT-GUE}}) &	\scriptsize({{\tt DIN-GUE}})&	\scriptsize({{\tt TAI-LAB}})&	
\scriptsize({{\tt WOO-LAB}}) &	
\scriptsize({{\tt ROM-REC}})&	
\scriptsize({{\tt ROM-TEC}})\\
 \hline\hline
$\mathbf{x}_1$	&	$\times$	&	$l_1t_1$	&	$l_1t_1$	&	$l_2t_3$	&	$l_1t_0$	&	$l_1t_0$	&	$l_1t_0$	&	$l_2t_0$	&	$l_2t_1$	&	$l_2t_1$	\\
$\mathbf{x}_2$	&	$\times$	&	$l_1t_1$	&	$l_1t_0$	&	$l_2t_3$	&	$l_1t_0$	&	$l_1t_1$	&	$l_1t_0$	&	$l_2t_0$	&	$l_1t_0$	&	$l_2t_1$	\\
$\mathbf{x}_3$	&	$l_1t_3$	&	$l_1t_1$	&	$l_1t_0$	&	$\times$	&	$l_2t_1$	&	$l_1t_1$	&	$l_1t_0$	&	$l_2t_0$	&	$l_2t_0$	&	$l_1t_0$	\\
$\mathbf{x}_4$	&	$l_1t_3$	&	$l_1t_1$	&	$l_1t_0$	&	$\times$	&	$l_2t_1$	&	$l_1t_1$	&	$l_1t_0$	&	$l_2t_0$	&	$l_1t_0$	&	$l_1t_0$	\\
$\mathbf{x}_5$	&	$\times$	&	$l_1t_1$	&	$l_1t_0$	&	$\times$	&	$l_1t_2$	&	$l_1t_2$	&	$l_1t_1$	&	$l_2t_0$	&	$l_1t_3$	&	$l_1t_0$	\\
$\mathbf{x}_6$	&	$\times$	&	$l_1t_3$	&	$l_1t_0$	&	$\times$	&	$l_1t_1$	&	$l_1t_2$	&	$l_1t_0$	&	$l_2t_0$	&	$l_2t_1$	&	$l_2t_2$	\\
$\mathbf{x}_7$	&	$\times$	&	$l_1t_1$	&	$l_1t_0$	&	$\times$	&	$l_1t_1$	&	$l_1t_2$	&	$l_1t_2$	&	$l_2t_0$	&	$l_2t_3$	&	$l_1t_0$	\\
$\mathbf{x}_8$	&	$\times$	&	$l_1t_1$	&	$l_1t_0$	&	$\times$	&	$l_2t_2$	&	$l_1t_3$	&	$l_1t_0$	&	$l_2t_0$	&	$l_1t_2$	&	$l_1t_1$	\\

 \end{tabular}
 \caption{Plans presented to the DM during the first iteration }\label{tab:firststrategies}
\end{table}

\subsection{Ranking of the proposed plans and elicitation of the preferences}

The DM, faced with the plans in Table \ref{tab:firststrategies}, pointed out that there are some priorities and requirements to bear in mind: 
\begin{itemize}
	\item The tailor's laboratory ({{\tt TAI-LAB}}), which also contains the laundry, must be built immediately so that the residents can be accommodated. This service cannot be outsourced because it is based on the crucial principles of Ecovillage, such as water recycling.
	\item In the identified plans, a mixed-use of kitchens and refectories for guests and residents is implemented at the starting period $t_0$: the DM considered this to be very reasonable. From a strategic point of view, the DM points out that it makes sense to have alternatives where guest kitchens are implemented initially since there may be catering without residents initially, but not vice versa.
\item Preference has to be given to plans where the recreational room ({{\tt ROM-REC}}) is in the Upper Borough, where all other functions are because it is more convenient for guests. In overnight accommodation, the spaces can be used interchangeably between residents and external guests, in the first phase of the settlement, there is a high degree of adaptability because guests and residents are not very dissimilar. Again, the above requirements were considered by adding corresponding constraints to the optimization problems to be solved to define the plans to propose to the DM.
\end{itemize}

More in detail, commenting on the first four plans related to the budget $B_1$, the DM observes that plans $\mathbf{x}_1$ is preferred over plan $\mathbf{x}_2$ because the kitchen ({{\tt KIT-GUE}}) and guest dining room ({{\tt DIN-GUE}}) are located in a building that is most suitable for hospitality in the medium to long term; plan $\mathbf{x}_4$ is preferred over plan $\mathbf{x}_3$ because the recreational room ({{\tt ROM-REC}}) is located in the upper borough, which is more convenient for short-stay guests. The DM also underlined that plan $\mathbf{x}_1$ is preferred to plan $\mathbf{x}_3$, considering the higher income it could provide since the catering can be realized here immediately. Then, applying the deck of the cards method, we asked the DM to rank order the  plans related to budget $B_1$, also providing a measure of the strength of preferences in terms of the number of blank cards between each plan and the following one in the preference ranking. The DM provided the following ranking, identified with $R_{50}$ with the number of blank cards shown between parenthesis $[\;]$, with ${{\tt{x_0^1}}}$  representing a fictitious plan identifying a zero level for budget $B_1$:

$$  {{\tt {x_0^1} }} \;\; [5] \;\; {{\tt x_3}} \;\; [0] \;\;{{\tt x_4}} \;\; [2] \;\;{{\tt x_2}} \;\; [3] \;\;{{\tt x_1}}$$

Commenting on the plans for the budget configuration $B_2$, the DM states that these ones are less preferred because there is no residence facility in any of them. Plan $\mathbf{x}_6$ is the preferred one because it has selected the kitchen for guests ({{\tt KIT-GUE}}) and a refectory ({{\tt DIN-GUE}}). For the guests, the most connotative room is the one for recreational activities ({{\tt ROM-REC}}), which are more rare and uncommon for the region (such as yoga, martial arts) and, together with the dining activity, are also the most profitable. The worst plan is $\mathbf{x}_8$ because it does not schedule the opening of the technical room ({{\tt ROM-TEC}}) at the starting period; plan $\mathbf{x}_7$ is worse than plan $\mathbf{x}_5$ because there is no tailoring laboratory ({ {\tt TAI-LAB}}).  Then, we asked the DM to rank order the plans and to insert blank cards representing the strength of preferences concerning plans related to budget $B_2$. The DM provided the following preference information, with ${{\tt {x_0^2}}}$ representing a fictitious plan identifying a zero level for budget $B_2$:
$$  {{\tt {x_0^2}}} \;\; [2] \;\; {{\tt x_8}} \;\; [3] \;\;{{\tt x_7}} \;\; [2] \;\;{{\tt x_5}} \;\; [5] \;\;{{\tt x_6}}$$ 
%\begin{table}[htb!]
%\centering
 %\begin{tabular}{cc|cc}
 %\multicolumn{2}{c|}{$F_{B_1}$}&\multicolumn{2}{c}{$F_{B_2}$}\\
 %\hline\hline
 %Plans&No. of cards&Plans&No. of cards\\
 %$x_6$&$x_1$&3\\
%$[5]$&[3]\\ 
%$x_5$&$x_2$\\
%$[2]$&[2]$\\ 
%$x_7$&3&$x_4$&0\\
 %$x_8$&2&$x_3$&5\\
%\end{tabular}
 %\caption{Partial ranking suggested by the DM during the first interaction }\label{tab:planspartial}
%\end{table}
To create a single ranking between the plans related to the budget configuration $B_1$ (considered in general favorite) and the plans related to the budget configuration $B_2$, we asked the DM to define a number of cards between the worst plan related to $B_1$, that is $\mathbf{x_3}$, and the best plan related to $B_2$, that is $x_6$. The DM established a distance of seven cards, justifying this  significant distance considering that plans related to budget configuration $B_2$ do not present any housing functions, which would mean creating more of a restaurant with some related activities than a real Ecovillage. In addition, if, on the one hand, the first four plans require twice the budget of the others, on the other hand, they provide more than double the revenue. The final ranking with the related cards measuring the strength of preferences between one plan and the following ones, identified with the following preference information $R_{Tot}$,  where ${{\tt{x_0}}} = {{\tt {x_0^2}}}$ is interpreted as a general zero level:
$$ {{\tt x_0}} \;\; [2]\;\; {{\tt x_8}}\;\; [3]\;\; {{\tt x_7}} \;\; [2] \;\;{{\tt x_5}} \;\; [5] \;\;{{\tt x_6}} \;\; [7]  \;\; {{\tt x_3}} \;\; [0] \;\; {{\tt x_4}} \;\; [2] \;\;{{\tt x_2}} \;\; [3] \;\;{{\tt x_1}} \;\; $$

%%\begin{figure}[htp!]
%% \centering
%%\includegraphics{Cards First.jpg}
%%%% \caption{Ranking of the first eight selected plans with blank cards measuring strength of preferences
%%}\label{fig:Mappa}

%%\end{figure}

Using the above preference information supplied by the DM in terms of ranking of plans, we induced the parameters of a more complex value function taking into account the interaction between criteria and the synergy between projects. More precisely, we proceeded as follows. We considered a value function $U(\mathbf{x})$ expressed in terms of a Choquet integral aggregating evaluations on the already considered four criteria $g_1, g_2, g_3$ and $g_4$ plus the further criterion $syn$ taking value 1 if in the considered plan there is the synergy between facilities and zero vice versa. The criterion $syn$ was added because the DM considered a specific relevance to the interaction between facilities ({{\tt KIT-GUE}}) and ({{\tt DIN-GUE}}) larger than the increase $\sigma_r$ given to the evaluation of considered facilities on considered criteria. We took into consideration the interaction between pairs of the four criteria $g_1, g_2, g_3$ and $g_4$, while we did not consider any interaction between the synergy $syn$ and one of the criteria $g_1, g_2, g_3$ and $g_4$. In consequence, the value function we adopted has the following formulation
$$U(\mathbf{x})=\sum_{j=1}^4 w_j g_j(\mathbf{x})+\sum_{j,j'=1,2,3,4, j\neq j'} w_{jj'}min(g_j(\mathbf{x}),g_{j'}(\mathbf{x}))+w_{syn}syn(\mathbf{x})$$ 
with $\sum_{j=1}^4 w_j+\sum_{j,j'=1,2,3,4, j\neq j'} w_{jj'}+w_{syn}=1$, $w_{syn} \geqslant 0$, and $w_j, j=1,2,3,4,$ and $w_{j,j'},j,j'=1,2,3,4, j\neq j'$, satisfying all the constraint of the Choquet non-additive weights.  
We applied the DOR methodology to the preference information given by the DM in terms of SRFII deck of the cards method to 
\begin{enumerate}
\item the ranking of plans related to budget $B_2$, identified as $R_{50}$;	
\item the ranking of plans related to budget $B_1$ identified as $R_{100}$,	
 	\item the whole ranking of plans related to budget $B_1$ and $B_2$, identified as $R_{Tot}$.
	\end{enumerate}
Then, formulating the problem in terms of the LP  (19) in Section 2, we computed three vectors of non-additive weights, reported in Table \ref{tab:weightssecond},  
%%$$\mathbf{\hat{w}^h}=[w^h_1,w^h_2,w^h_3,w^h_4,w^h_{1,2},w^h_{1,3}w^h_{1,4}w^h_{2,3}w^h_{2,4}w^h_{3,4},w^h_{syn}], h=1,2,3$$
 for the Choquet integral formulation of the value function $U(\mathbf{x})$ corresponding to the ranking obtained with the deck of cards method.

\begin{table}[htb!]
\centering
 \begin{tabular}{c|ccccccccccc}
&$w_1$&	$w_2$&	$w_3$&	$w_4$&	$w_{12}$&	$w_{13}$&	$w_{14}$&	$w_{23}$&	$w_{24}$&	$w_{34}$&	$w_{sin}$\\

 \hline\hline
$\bold{{w}^{R_{50}}}$&	0.05&	0&	0.502&	0&	0&	0&	0&	0&	0.175&	0&	0.273\\
$\bold{{w}^{R_{100}}}$&	0&	0&	0&	0&	0&	0&	0.468&	0&	0&	0&	0.532\\
$\bold{{w}^{R{Tot}}}$&	0.306&	0&	0.455&	0&	0&	0&	0&	0&	0&	0&	0.239\\
 \end{tabular}
 \caption{Nonadditive weights for the value function expressed in terms of a Choquet integral}\label{tab:weightssecond}
\end{table}

In Tables \ref{tab:regressionmodel_Choquet_PrimaIterazione_R50}, \ref{tab:regressionmodel_Choquet_PrimaIterazione_R100} and \ref{tab:regressionmodel_Choquet_PrimaIterazione_RTot} we reported the values assigned to each plan with the deck of cards method, the value function $U$, the  corrected value function and the differences $\sigma^+$ and $\sigma^-$ for each of the configuration introduced, respectively.

\begin{table}[!h]
\begin{center}
\caption{Scores assigned to plans by the value function $U$ obtained solving the LP problem  (19) for ranking $R_{50}$}\label{tab:regressionmodel_Choquet_PrimaIterazione_R50}
		\begin{tabular}{cccccc}
	\hline
   Plans  & $U(\mathbf{x}_i)$ & $\nu(\mathbf{x}_i)$ & $k\cdot\nu(\mathbf{x}_i)$& $\sigma^+(\mathbf{x}_i)$ & $\sigma^-(\mathbf{x}_i)$ \\
	\hline
${ {\tt x_5}}$ &0.31&	10&	0.31&	0&	0\\

   ${ {\tt x_6}}$&0.5&	16&	0.5&	0&	0\\

${ {\tt x_7}}$& 0.22&	7&	0.22&	0&	0\\

${ {\tt x_8}}$ &0.09&	3&	0.09&	0&	0\\

		\hline
		\end{tabular}
\end{center}
\end{table}
\begin{table}[!h]
\begin{center}
\caption{Scores assigned to plans by the value function $U$ obtained solving the LP problem  (19) for ranking $R_{100}$}\label{tab:regressionmodel_Choquet_PrimaIterazione_R100}
		\begin{tabular}{cccccc}
	\hline
   Plans  & $U(\mathbf{x}_i)$ & $\nu(\mathbf{x}_i)$ & $k\cdot\nu(\mathbf{x}_i)$& $\sigma^+(\mathbf{x}_i)$ & $\sigma^-(\mathbf{x}_i)$ \\
	\hline
${ {\tt x_1}}$ &0.53&	14&	0.53&	0&	0\\
  
   ${ {\tt x_2}}$&0.70&	10&	0.38&	0&	0.32\\

${ {\tt x_3}}$& 0.25&	6&	0.23&	0&	0.02\\

${ {\tt x_4}}$ &0.27&	7&	0.27&	0&	0\\

		\hline
		\end{tabular}
\end{center}
\end{table}

\begin{table}[!h]
\begin{center}
\caption{Scores assigned to plans by the value function $U$ obtained solving the LP problem  (19) for ranking $R_{Tot}$}\label{tab:regressionmodel_Choquet_PrimaIterazione_RTot}
		\begin{tabular}{cccccc}
	\hline
   Plans  & $U(\mathbf{x}_i)$ & $\nu(\mathbf{x}_i)$ & $k\cdot\nu(\mathbf{x}_i)$& $\sigma^+(\mathbf{x}_i)$ & $\sigma^-(\mathbf{x}_i)$ \\
	\hline
${ {\tt x_1}}$ &0.89&	32&	0.89&	0&	0\\
${ {\tt x_2}}$&0.96&	28&	0.78&	0&	0.18\\
${ {\tt x_3}}$& 0.72&	24&	0.67&	0&	0.06\\
${ {\tt x_4}}$ &0.7&	25&	0.7&	0&	0\\

${ {\tt x_5}}$ &0.28&	10&	0.28&	0&	0\\
${ {\tt x_6}}$ &0.11&	16&	0.45&	0.34&	0\\
${ {\tt x_7}}$ &0.06&	7&	0.19&	0.14&	0\\
${ {\tt x_8}}$ &0.06&	3&	0.08&	0.03&	0\\

		\hline
		\end{tabular}
\end{center}
\end{table}
\subsection{Definition of a new set of plans}
At this point, on the basis of the discussion conducted with the DM we generate a new set of plans optimizing a value function $U$ formulated in terms of a Choquet integral  related to the weight vectors $\mathbf{w^{R_{50}}},\mathbf{w^{R_{100}}}$ and $\mathbf{w^{R_{Tot}}}$ induced in the previous step. We consider the two budget configurations $B_1$ and $B_2$ as previously defined. We also impose the constraints that at least a kitchen should be selected and that the facility ({{\tt TAI-LAB}}) should be selected earlier than the facilities ({{\tt RES-WWO}}) and  ({{\tt WOO-LAB}}), according to the preferences expressed by the DM during this second discussion. 
 We also included a plan for each of the budget configurations with the complete order and with an additional constraint on the presence of at least  one of the residences; this is to investigate if the DM would prefer such plans that will allow him since the beginning to host some guests in the Ecovillage. The synergy constraint, related to the activation of facilities ({{\tt KIT-GUE}}) and ({{\tt DIN-GUE}}), was always included according to the DM preferences expressed in the previous step. In total, we generated eight plans, obtained by the combinations of the two budget scenarios, the three sets of weights $\bold{w^{R50}}$, $\bold{w^{R100}}$ and $\bold{w^{R_{Tot}}}$ and the presence of at least one of the residences with the set of weights $\bold{w^{R_{Tot}}}$. 
More precisely, the  plans selected could be obtained in this way
 
 \begin{itemize}
     \item $\mathbf{x}_1^{'}$, for budget $B_1$, weight vector $\mathbf{w^{R_{50}}}$;
      \item $\mathbf{x}_2^{'}$, for budget $B_1$, weight vector $\mathbf{w^{R_{100}}}$;
      \item $\mathbf{x}_3^{'}$, for budget $B_1$, weight vector $\mathbf{w^{R_{Tot}}}$;
      \item $\mathbf{x}_4^{'}$, for budget $B_1$, weight vector $\mathbf{w^{R_{Tot}}}$, with the residence constraint;
      \item $\mathbf{x}_5^{'}$, for budget $B_2$, weight vector $\mathbf{w^{R_{100}}}$;
      \item $\mathbf{x}_6^{'}$, for budget $B_2$, weight vector $\mathbf{w^{R_{100}}}$ ;
      \item $\mathbf{x}_7^{'}$, for budget $B_2$, weight vector $\mathbf{w^{R_{Tot}}}$;
      \item $\mathbf{x}_8^{'}$, for budget $B_2$, weight vector $\mathbf{w^{R_{Tot}}}$, with the residence constraint.

 \end{itemize}
 
 Those new plans were then presented to the DM and reported in Table \ref{tab:secondstrategies}.

\begin{table}[htb!]
\centering
 \begin{tabular}{c|cccccccccc}

  & \scriptsize({ {\tt RES-WWO}}) & \scriptsize({{\tt KIT-WWO}})&\scriptsize({{\tt REF-WWO}})&\scriptsize({{\tt ROM-GUE}}) &
\scriptsize({{\tt KIT-GUE}}) &	\scriptsize({{\tt DIN-GUE}})&	\scriptsize({{\tt TAI-LAB}})&	
\scriptsize({{\tt WOO-LAB}}) &	
\scriptsize({{\tt ROM-REC}})&	
\scriptsize({{\tt ROM-TEC}})\\
 \hline\hline
$\mathbf{x}_1^{'}$	&$\times$&	$l_1t_1$&	$l_1t_1$&	$l_1t_3$&	$l_1t_0$&	$l_1t_0$&	$l_1t_0$&	$l_1t_3$&	$l_2t_1$&	$l_1t_0$	\\
$\mathbf{x}_2^{'}$&	$\times$&	$l_1t_0$&	$l_1t_0$&	$l_1t_3$&	$l_2t_1$&	$l_1t_1$&	$l_1t_1$&	$l_1t_3$&	$l_1t_0$&	$l_1t_0$	\\
$\mathbf{x}_3^{'}$	&$\times$&	$l_1t_1$&	$l_1t_1$&	$l_2t_3$&	$l_1t_0$&	$l_1t_0$&	$l_1t_0$&	$l_1t_3$&	$l_2t_1$&	$l_1t_0$	\\
$\mathbf{x}_4^{'}$	&$\times$&	$l_1t_1$&	$l_1t_1$&	$l_2t_3$&	$l_1t_0$&	$l_1t_0$&	$l_1t_0$&	$l_1t_3$&	$l_2t_1$&	$l_1t_0$	\\
$\mathbf{x}_5^{'}$	&$\times$&	$l_1t_2$&	$l_1t_0$&	$\times$&	$l_1t_0$&	$l_1t_1$&	$l_1t_1$&	$l_1t_3$&	$l_2t_2$&	$l_2t_3$	\\
$\mathbf{x}_6^{'}$	&$\times$&	$l_1t_2$&	$l_1t_0$&	$\times$&	$l_1t_0$&	$l_1t_1$&	$l_1t_1$&	$l_1t_3$&	$l_2t_2$&	$l_2t_3$	\\
$\mathbf{x}_7^{'}$	$\times$&	$l_1t_3$&	$l_1t_1$&	$\times$&	$l_1t_0$&	$l_1t_2$&	$l_1t_0$&	$l_1t_3$&	$l_1t_1$&	$l_1t_1$	\\
$\mathbf{x}_8^{'}$	&$\times$&	$\times$&	$\times$&	$l_1t_3$&	$l_1t_0$&	$\times$&	$\times$&	$\times$&	$\times$&	$\times$	\\

 \end{tabular}
 \caption{Plans presented to the DM during the second iteration }\label{tab:secondstrategies}
\end{table}

The DM expressed his preference for the plan $\mathbf{x}_1^{'}$. He points out that the only inconsistency is that the recreation room ({{\tt ROM-TEC}}) in the new pavilion is too far away. 

In this sense, the DM detailed that the recreational room should ({{\tt ROM-REC}}) be close to the guest refectory ({{\tt DIN-GUE}}) (which, in turn, must be close to the guest kitchen ({{\tt KIT-GUE}})), and that the space must not be less than 30 m². %%Questo vincolo non dipende da dove sono definite le locations?
Otherwise, everything is congruent and the principle of environmental protection is respected.
With respect to the plan obtained with budget configuration $B_2$, the DM underlines that even considering the actual economic difficulties in starting the transformation process of the area, it constitutes a ``horizontal cut” that implies no overnight hospitality solution: having only the function ({{\tt DIN-GUE}}) is not interesting enough. More generally,  the DM expressed a preference of having at least two functions must be located in each of the transformed buildings. Therefore, we formulated these constraints and adopted the same weight vector $w^{R_{Tot}}$ for the budget configuration $B_1$ and adoption $\mathbf{w}^{R_{Tot}}$ that produced the preferred plan for the DM in the previous step, that is, $\mathbf{x'_1}$. A total of three new plans was generated, as follows: 
\begin{itemize}
    \item plan $\mathbf{x}_1^{''}$, obtained imposing that functions ({{\tt WOO-LAB}}) and ({{\tt ROM-REC}}) should not be both located in Location 1;
    \item plan $\mathbf{x}_2^{''}$, obtained imposing that in each building in which a function is activated, at least two functions are activated ; 
    \item plan $\mathbf{x}_3^{''}$,  obtained, imposing that at least two functions must be activated in each building.
\end{itemize}

\begin{figure}[htp!]
 \centering
\includegraphics{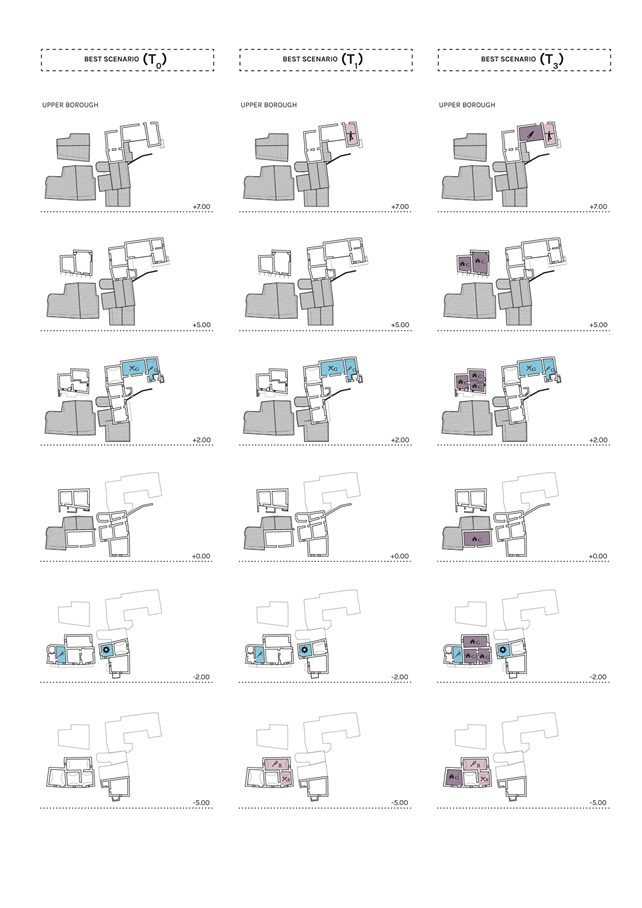}
 \caption{Selected facilities and their timing for the most representative plan $\mathbf{x}_2^{''}$
}\label{fig:Mappafinale}

\end{figure}

\begin{table}[htb!]
\centering
 \begin{tabular}{c|cccccccccc} 
  & \scriptsize({ {\tt RES-WWO}}) & \scriptsize({{\tt KIT-WWO}})&\scriptsize({{\tt REF-WWO}})&\scriptsize({{\tt ROM-GUE}}) &
\scriptsize({{\tt KIT-GUE}}) &	\scriptsize({{\tt DIN-GUE}})&	\scriptsize({{\tt TAI-LAB}})&	
\scriptsize({{\tt WOO-LAB}}) &	
\scriptsize({{\tt ROM-REC}})&	
\scriptsize({{\tt ROM-TEC}})\\
 \hline\hline
$\mathbf{x}_1^{''}$	&$\times$	&	$l_1t_1$	&	$l_1t_1$	&	$l_2t_3$	&	$l_1t_0$	&	$l_1t_0$	&	$l_1t_0$	&	$\times$	&	$l_2t_1$	&	$l_1t_0$	\\
$\mathbf{x}_2^{''}$&	$\times$	&	$l_1t_1$	&	$l_1t_1$	&	$l_2t_3$	&	$l_1t_0$	&	$l_1t_0$	&	$l_1t_0$	&	$l_1t_3$	&	$l_1t_1$	&	$l_1t_0$	\\
$\mathbf{x}_3^{''}$	&$\times$	&	$l_1t_3$	&	$l_1t_1$	&	$\times$	&	$l_1t_0$	&	$l_1t_2$	&	$l_1t_2$	&	$l_1t_3$	&	$l_2t_1$	&	$l_2t_0$	\\
 \end{tabular}
 \caption{Strategies presented to the DM during the third iteration }\label{tab:thirdstrategies}
\end{table}
Observing plan $\mathbf{x}_1^{''}$, the DM noted that there is a compact timing of the renovations while the locations are acceptable. He also pointed out that  there are only two critical points: the recreational room ({{\tt ROM-REC}})  remains disconnected from the transformed village and there is no woodworking room ({{\tt WOO-LAB}}).
Plan $\mathbf{x}_2^{''}$ is the most interesting for the DM because it locates everything in the borough above, it is compact and makes it comfortable to manage the space for guests and residents, and it has all the functions. There is a problem that the woodworking room ({{\tt WOO-LAB}}) is too close to the recreational room ({{\tt ROM-REC}}), so this location should be changed. %%siamo sicuri che dobbiamo/vogliamo inserire questo commento? 
Plan $\mathbf{x}_3^{''}$ is the least preferred, especially concerning the timing of the implementation of the various functions, with some functions having to be done together (e.g., there is the food serving space away from the kitchens). Therefore, the DM selected a plan 
$\mathbf{x}_2^{''}$ as the most representative one of his opinions. Let us also note that we interacted with the DM thanks to the use of technical representation of the Ecovillage in which the selected facilities and their timing were represented. For example, in Figure  \ref{fig:Mappafinale}, we report a representation of the selected facilities for the most representative plan for the DM.
To facilitate the DM's understanding of the temporal sequence of the realization of the facilities, different colors have been used: facilities activated at $t_1$ are in light blue, those at $t_2$ in pink and those at $t_3$ are lilac. The arrangement of the floor plans made the communication and evaluation of the different plans particularly effective.

%\begin{figure}[htp!]\label{fig:finalesoluzione}
 %\centering
%\includegraphics{SoluzioneFinale.png}
 %\caption{The preferred strategy with the functions distributed at the respective times and locations
%}
 %\label{fig:finalesoluzione}
%\end{figure}

\section{Conclusions}\label{Conc}

In this paper, we proposed a new procedure to support territorial planning decisions formulated in terms of multiobjective optimization problems that, in the framework of the so-called space-time model, define which facilities to implement in which locations and at which time \citep{barbati2020general} A relevant point of the proposed procedure is the consideration of the interaction between criteria through introducing the Choquet integral in the space-time model. However, taking into account the relevance of interaction with decision makers, the main point of the proposed methodology is a new elicitation procedure called Deck of the cards based Ordinal Regression (DOR) that put together the easy understandability for the DM of  the deck of the cards method and the versatility of ordinal regression to induce preference parameters of a specific value function representing DM's preferences. In the proposed decision-aiding procedure to support territorial planning, DOR method is used to construct a value function representing DM's preferences when dealing with multiobjective optimization problems. However, the DOR method is general and can be applied to any decision problem. 

We illustrated the application of the whole procedure in a case study related to the foundation of an Ecovillage. In this specific context, the challenge is to create an environmentally responsible settlement, able to reconcile two antithetical values such as the desire to pursue an informal economy that is entirely unrelated to commercial logic and, at the same time the need to achieve economic self-sufficiency of the settlement. In addition, there is the presence of three types of users, i.e. the residents, the WWOOFERs, and the guests, imposing location choices with very different timeframes (short, medium and long term), relating both to the construction of the various buildings and to the subsequent management of the functions to be carried out in them. This type of application is specifically relevant because it can be viewed as a case study for decision-aiding related to choices involving aspects such as sustainability and social responsibility which are more and more fundamental for future generations' planet Earth.  
With respect to future developments of the research, the following points seem to us the most promising:
\begin{itemize}
\item The territorial decision support methodology we are proposing could be applied in other contexts, and  different decision-aiding problems could be considered.
\item The proposed methodology could be integrated to include the opinions of several DMs and could be adapted in a group context decision-making process.
\item  Applications of the methodology to large-scale planning could be developed.
\item Theoretical advances to consider much longer time periods, concerning also intergenerational
issues could be dealt with.
\end{itemize}
Finally, we want to point out that a specific interest is related to the DOR methodology that could be tested in several diversified decision problems to verify its advantages in real-life applications.

\newpage
\bibliographystyle{model2-names}
\bibliography{Sample.bib}
\addcontentsline{toc}{section}{References}

\end{document}